\documentclass{amsart}

\usepackage[
	backend=biber,
	style=alphabetic,
	backref=true,
	url=false,
	doi=false,
	isbn=false,
	eprint=false]{biblatex}

\renewbibmacro{in:}{} 
\AtEveryBibitem{\clearfield{pages}} 

\DeclareFieldFormat{title}{\myhref{\mkbibemph{#1}}}
\DeclareFieldFormat
[article,inbook,incollection,inproceedings,patent,thesis,unpublished]
{title}{\myhref{\mkbibquote{#1\isdot}}}

\newcommand{\doiorurl}{%
	\iffieldundef{url}
	{\iffieldundef{eprint}
		{}
		{http://arxiv.org/abs/\strfield{eprint}}}
	{\strfield{url}}%
}

\newcommand{\myhref}[1]{%
	\ifboolexpr{%
		test {\ifhyperref}
		and
		not test {\iftoggle{bbx:eprint}}
		and
		not test {\iftoggle{bbx:url}}
	}
	{\href{\doiorurl}{#1}}
	{#1}%
}

\usepackage[
	bookmarks=true,
	linktocpage=true,
	bookmarksnumbered=true,
	breaklinks=true,
	pdfstartview=FitH,
	hyperfigures=false,
	plainpages=false,
	naturalnames=true,
	colorlinks=true,
	pagebackref=false,
	pdfpagelabels]{hyperref}

\hypersetup{
	colorlinks,
	citecolor=NavyBlue,
	filecolor=NavyBlue,
	linkcolor=NavyBlue,
	urlcolor=NavyBlue
}

\usepackage[capitalize,noabbrev]{cleveref}
\crefname{subsection}{\S\kern-.3em}{subsections}
\Crefname{subsection}{\S\kern-.3em}{Subsections}

\calclayout

\makeatletter
\@namedef{subjclassname@2020}{%
	\textup{2020} Mathematics Subject Classification}
\makeatother

\usepackage[T1]{fontenc}
\usepackage{lmodern}
\usepackage{microtype}

\usepackage{mathtools}
\usepackage{amssymb}
\usepackage{mathbbol} 

\usepackage[dvipsnames]{xcolor}

\usepackage{tikz-cd}

\usepackage{csquotes}
\usepackage[shortcuts]{extdash} 

\usepackage{enumerate}
\usepackage{enumitem}

\usepackage{todo}

\newtheorem{theorem}[subsection]{Theorem}
\newtheorem*{theorem*}{Theorem}

\newtheorem*{proposition*}{Proposition}
\newtheorem{lemma}[subsection]{Lemma}
\newtheorem*{lemma*}{Lemma}
\newtheorem{corollary}[subsection]{Corollary}
\newtheorem*{corollary*}{Corollary}

\theoremstyle{definition}

\newtheorem*{definition*}{Definition}
\newtheorem{remark}[subsection]{Remark}
\newtheorem*{remark*}{Remark}

\newtheorem*{example*}{Example}

\newtheorem*{construction*}{Construction}

\newtheorem*{convention*}{Convention}

\newtheorem*{terminology*}{Terminology}

\newtheorem*{notation*}{Notation}

\newtheorem*{question*}{Question}

\DeclareMathOperator{\bd}{\partial}

\newcommand{\ot}{\otimes}
\DeclareMathOperator{\EZ}{EZ}
\DeclareMathOperator{\AW}{AW}

\newcommand{\N}{\mathbb{N}}

\newcommand{\Sym}{\mathbb{S}}

\newcommand{\Ftwo}{{\mathbb{F}_2}}

\newcommand{\simplex}{\triangle}

\DeclareMathOperator{\Sq}{Sq}
\DeclareMathOperator{\ind}{ind}

\DeclareMathOperator{\chains}{N}
\DeclareMathOperator{\cochains}{N^{\vee}}

\DeclarePairedDelimiter\bars{\lvert}{\rvert}

\DeclarePairedDelimiter\set{\{}{\}}

\newcommand{\id}{\mathsf{id}}

\newcommand{\Hom}{\mathrm{Hom}}

\colorlet{DefinitionGreen}{PineGreen!100!black}
\newcommand{\defn}[1]{{\color{DefinitionGreen}\emph{#1}}}

\newcommand{\rH}{\mathrm{H}}

\newcommand{\rP}{\mathrm{P}}

\newcommand{\rR}{\mathrm{R}}
\newcommand{\rS}{\mathrm{S}}

\newcommand{\cE}{\mathcal{E}}

\newcommand{\cP}{\mathcal{P}}

\newcommand{\cX}{\mathcal{X}}

\newcommand{\canonical}{\triangle^{\!\mathrm{cn}}}
\DeclareMathOperator{\pired}{\pi_{\text{red}}}
\DeclareMathOperator{\piirred}{\pi^{\perp}_{\text{red}}}
\newcommand{\F}{\Ftwo}
\renewcommand{\P}{\mathrm{P}}

\DeclareMathOperator*{\displaytensor}{\otimes}
\let\union\cup
\renewcommand{\cup}{\smallsmile}

\renewcommand{\simplex}{\mathbb{\Delta}}
\newcommand{\barxi}{{\bar{\xi}}}

\DeclareMathOperator{\TR}{TR}
\newcommand{\coeff}[2]{\langle #1 ; #2 \rangle}
\newcommand{\AWd}{{\mathrm{AW}}}
\title[Axiomatic characterization of cup-$i$ products]{An axiomatic characterization of Steenrod's cup-$i$ products}

\author{Anibal~M.~Medina-Mardones}
\address{Department of Mathematics, Western University, Canada}
\email{\href{mailto:anibal.medina.mardones@uwo.ca}{anibal.medina.mardones@uwo.ca}}

\date{\today}

\begin{document}

\begin{abstract}
	We show that any natural construction of cup-$i$ products on the normalized cochains of simplicial sets, parameterized by the minimal free resolution of the trivial \(\Ftwo[\Sym_2]\)-module, is isomorphic--not just homotopic--to Steenrod's original construction if it is non-zero, irreducible, and free.
	We show that these three conditions are independent, and use them to prove that all \mbox{cup-$i$} constructions in the literature represent the same isomorphism class.
\end{abstract}

	\maketitle
	\tableofcontents

\section{Introduction}\label{s:introduction}

In his seminal paper \cite{steenrod1947products}, Steenrod used explicit formulas to introduce on the cochains of triangulated spaces the \textit{cup-$i$ products}
\[
\cup_i\ \colon \cochains(X; \F)^{\ot 2} \to \cochains(X; \F).
\]
These bilinear maps give rise to the celebrated \textit{Steenrod squares}
\[
\begin{tikzcd}[column sep=0, row sep=0]
	\Sq^k \colon &[-5pt] \rH^\vee(X; \F) \rar & \rH^\vee(X; \F) \\
	& {[\alpha]} \rar[maps to] & \big[(\alpha \cup_{-\bars{\alpha}-k} \alpha)\big]
\end{tikzcd}
\]
lying at the heart of stable homotopy theory.

Steenrod's $\cup_0$ product agrees with the Alexander--Whitney product on cochains, inducing the commutative product structure on cohomology.
The higher cup-$i$ products can be interpreted as coherent homotopies witnessing the derived commutativity of $\cup_0$ at the cochain level.

Later work by Steenrod and others established the existence of the cup-$i$ products and, consequently, of the Steenrod squares using an indirect argument based on acyclic carrier methods.
This approach became standard since any set of choices for the cup-$i$ products homotopic to Steenrod's original gives rise to the same cohomology operations, which had been axiomatically characterized.
Consequently, the interest in specific formulas for the cup-$i$ products diminished in the second half of the century.

Recent interest in explicit cochain-level constructions has been driven in part by applications to topological phases of matter \cite{chen2012symmetry, chen2013symmetry, gu2014symmetry, kapustin2015fermionic_cobordisms, gaiotto2016spin, kapustin2017fermionic, meng2018classification, wang2020construction, barkeshli2022classification, feng2026anyonic}.
These applications require not only cochain representatives of Steenrod operations, such as the cup-$i$ products, but also cochain-level witnesses for the Cartan and Adem relations \cite{medina2020cartan, medina2021adem, medina2025odd_cartan}.
Mathematically interpreted, such constructions are used to build explicit cochain approximations to the Postnikov towers of generalized cohomology theories (\cite{brumfiel2016pontrjagin, brumfiel2018pontrjagin}).

In this work, we give an axiomatic characterization of Steenrod's original cup-$i$ construction up to isomorphism, rather than merely up to homotopy.
We work with natural constructions parameterized by the minimal free resolution and impose three additional conditions: nontriviality, irreducibility with respect to subsimplices, and maximal freeness under transposition.
We prove that these three conditions are independent and that the cup-$i$ formulas appearing in \cite{steenrod1947products, real1996computability, gonzalez-diaz1999steenrod, mcclure2003multivariable, berger2004combinatorial, medina2020prop1, medina2025fast_sq} all represent the same isomorphism class.
Thus, the classical uniqueness of the induced cohomology operations is strengthened to a uniqueness statement at the level of explicit cochain constructions.

\subsection*{Related work}

The corresponding picture at odd primes is less complete.
Explicit cup-$(p,i)$ products were introduced in \cite{medina2021may_st} and implemented in the computer algebra system \href{https://comch.readthedocs.io/en/latest/}{\texttt{ComCH}} \cite{medina2021comch}.
A second family was introduced in \cite{medina2024connected}; for $p > 3$, it disagrees with the earlier construction.
For $p = 3$ no axiomatic characterization exists but no examples are known that rule out this possibility.

The present paper concerns simplicial cup-$i$ products.
Analogous axiomatic questions can be posed for cubes or more general products of simplices, for which explicit formulas are already known: see \cite{pilarczyk2016cubical, medina2022cube_einfty} for the cubical case and \cite{medina2024multisimplicial} for the multisimplicial one.


\subsection*{Outline}

We begin by presenting the notion of \mbox{cup-$i$} construction and our axiomatic characterization in \cref{s:statement}.
We postpone the proof of our main theorem until after \cref{s:reformulation}, where the statement is recast in more category theoretic terms.
\cref{s:irreducible} studies consequences of irreducibility, including that every irreducible semi-simplicial \mbox{cup-$i$} construction is simplicial, and \cref{s:freeness} records the key consequence of freeness.
\cref{s:proof} is devoted to the proof of our axiomatic characterization.
In \cref{s:examples} we show that the axioms of our characterization are independent.
In \cref{s:others} we show that all formulas in the literature give rise to isomorphic simplicial \mbox{cup-$i$} constructions.


\subsection*{Acknowledgment}

We would like to thank the anonymous referee as well as Mark Behrens, Greg Brumfiel, Tim Campion, Kathryn Hess, Federico Cantero Mor\'an, Riley Levy, John Morgan, Stephan Stolz, and Dennis Sullivan for their insightful suggestions and comments about this project.

The author is also grateful for the hospitality of the Topology and Neuroscience Laboratory at EPFL, the Max Planck Institute for Mathematics in Bonn, and LAGA at Universit\'e Sorbonne Paris Nord.

Partial support for this project was provided by grant ANR 20 CE40 0016 01 PROJET HighAGT and NSERC grants RES000678 and R7444A03

\section{Cup-$i$ constructions}\label{s:statement}

Let $\Ftwo$ be the field with two elements.
We denote the (normalized) cochains with $\Ftwo$ coefficients of a simplicial or semi-simplicial set $X$ by $\cochains(X)$, graded homologically and concentrated in non-positive degrees.
We remind the reader that semi-simplicial sets are simplicial sets without degeneracy maps or, more precisely, contravariant functors from the injective simplex category to the category of sets.

\subsection{Cup-$i$ products}

Let $\Sym_2$ be the group with one non-identity element $T$.
Let
\[
\begin{tikzcd}[column sep=20pt]
W = \Big(
\Ftwo[\Sym_2]\{e_0\} &
\Ftwo[\Sym_2]\{e_1\} \arrow[l, "\ 1+T"'] &
\arrow[l, "\ 1+T"'] \dotsb \Big)
\end{tikzcd}
\]
be the minimal free resolution of $\Ftwo$ by $\Ftwo[\Sym_2]$-modules.
We notice that the group of automorphisms of $W$ is isomorphic to $\prod_{\N} \Sym_2$ with any such $\phi$ determined by a choice, for every $i \in \N$, between $\phi(e_i) = e_i$ or $\phi(e_i) = Te_i$.

A \defn{\mbox{cup-$i$} product structure} on a chain complex $A$ is a chain map
\[
W \ot_{\F[\Sym_2]} A^{\ot 2} \to A
\]
where $T$ acts by transposition on $A^{\ot 2}$ and by right multiplication on $W$.
We denote the image of the class of $e_i \ot \alpha \ot \beta$ by $\alpha \cup_i \beta$.
Unpacking this structure we have the following defining identity for any $i \in \N$ and $\alpha, \beta \in A$:
\[
\alpha \cup_{i-1} \beta + \beta \cup_{i-1} \alpha =
\bd (\alpha \cup_{i} \beta) + (\bd\alpha) \cup_{i} \beta + \alpha \cup_{i} (\bd\beta)
\]
with the convention $\alpha \cup_{-1} \beta = 0$.
By setting $i = 0$ we see that $\cup_0$ is compatible with the differential.
The identity obtained by setting $i = 1$ implies that the product induced by $\cup_0$ in the homology of $A$ is commutative, with $\cup_1$ enforcing this relation at the chain level.
The degree 1 product need not itself be commutative, but its failure to be so is controlled by a homotopy defined by $\cup_2$, and so on.

An \defn{isomorphism} of \mbox{cup-$i$} product structures on $A$, or simply a \defn{cup-$i$ isomorphism}, is an automorphism $\phi$ of $W$ making the following diagram commute:
\begin{center}
	\begin{tikzcd}[column sep=5, row sep=15]
		W \displaytensor_{\F[\Sym_2]} A^{\ot 2} \arrow[dr, in=180, out=-90] \arrow[rr, "\phi \, \ot \, \id \, "] & &
		W \displaytensor_{\F[\Sym_2]} A^{\ot 2}. \arrow[dl, in=0, out=-90] \\
		& A &
	\end{tikzcd}
\end{center}

A \defn{(semi-)simplicial \mbox{cup-$i$} construction} is a \mbox{cup-$i$} product structure on $\cochains(X)$ for every (semi-)simplicial set $X$ that is natural with respect to (semi-)simplicial maps.
Similarly, an \defn{isomorphism} of (semi-)simplicial \mbox{cup-$i$} constructions is a natural family of cup-$i$ isomorphisms for each such $X$.
Notice that a simplicial \mbox{cup-$i$} construction is a semi-simplicial one satisfying additional constraints, i.e., compatibility with degeneracy maps.

\subsection{Axioms}

The first axiom alluded to in the abstract, naturality, has been explicitly included in our definition of (semi-)simplicial \mbox{cup-$i$} construction; whereas the second, minimality, is manifested in the use of $W$ instead of an arbitrary resolution of $\F$ as \(\F[\Sym_2]\)-module.

The \defn{zero \mbox{cup-$i$} construction} assigns the zero map
\[
W \ot_{\F[\Sym_2]} \cochains(X)^{\ot 2} \to \cochains(X)
\]
to every (semi-)simplicial set $X$.
Any other construction is said to be \defn{non-zero}.



Let $X$ be a semi-simplicial set.
Below we identify simplices with their associated chain and cochain basis elements.
A semi-simplicial \mbox{cup-$i$} construction is said to be \defn{irreducible} if for any simplex $x$ in $X$
\[
\boxed{\Big( y^{(1)} \cup_{i} y^{(2)} \Big)(x) = 0}
\]
whenever $y$ is a proper face of $x$ and $y^{(1)}$ and $y^{(2)}$ are faces of $y$.
It is \defn{free} if for any two simplices $x$ and $x'$ in $X$
\[
\boxed{x \cup_{i} x' = x' \cup_{i} x} \
\Longrightarrow \
\boxed{x \cup_{i} x' = 0}
\]
whenever $\bars{x} \neq i$ or $\bars{x'} \neq i$.

\subsection{Main Theorem}\label{ss:main_theorem}

\textit{There is up to isomorphism only one non-zero, irreducible and free semi-simplicial \mbox{cup-$i$} construction.
Furthermore, it is simplicial.}

\section{Reformulation}\label{s:reformulation}

This section provides a more categorical perspective on \mbox{cup-$i$} constructions, one that is also dual in two complementary senses:
linearly, we pass from cochains to chains, their pre-dual; and combinatorially, we index the faces of a simplex by their complementary vertex sets, as in Alexander duality.

\subsection{Cup-\textit{i} coproducts}

A \defn{cup-$i$ coproduct structure} on a chain complex $C$ is an $\Sym_2$-equivariant chain map
\[
\begin{tikzcd}[column sep=small,row sep=0]
	\triangle \colon &[-15pt] W \rar & \Hom(C, C^{\ot 2}) \\
	& e_i \rar[maps to] & \triangle_i.
\end{tikzcd}
\]
Such a map $\triangle$ is equivalent to a collection of linear maps $\set{\triangle_i}_{i\in\N}$ satisfying
\[
\bd \circ \, \triangle_i + \triangle_i \circ \bd =
(1+T) \triangle_{i-1}
\]
for all $i \in \N$ with the convention $\triangle_{-1} = 0$.

Let us denote $\Hom(C, \F)$ by $C^\vee$ and assume $C$ is finite-dimensional.
Using the hom-tensor adjunction and the finite dimensionality of $C$ we have
\begin{equation}\label{eq:hom-tensor}
	\begin{split}
		\Hom \big(W \ot_{\F[\Sym_2]} (C^\vee)^{\ot 2}, C^\vee \big) & \cong
		\Hom_{\F[\Sym_2]} \big( W, \Hom((C^\vee)^{\ot 2}, C^\vee) \big) \\ & \cong
		\Hom_{\F[\Sym_2]} \big( W, \Hom(C, C^{\ot 2}) \big)
	\end{split}
\end{equation}
as chain complexes of $\F$-modules.
In other words, the linear duality functor induces a bijection between \mbox{cup-$i$} product structures on $C^\vee$ and cup-$i$ coproduct structures on $C$.

%

\subsection{Representables}\label{ss:naturality}

A semi-simplicial \mbox{cup-$i$} construction $\triangle$ is completely determined by its restriction to representable semi-simplicial sets
\[
W \ot_{\Ftwo[\Sym_2]} \cochains(\simplex^n)^{\ot 2} \to \cochains(\simplex^n).
\]
Given that each of the chain complexes $\chains(\simplex^n)$ is finite-dimensional, $\triangle$ is determined by the natural collection of cup-$i$ coproduct structures
\[
W \to \Hom(\chains(\simplex^n), \chains(\simplex^n)^{\ot 2})
\]
defined by the hom-tensor adjunction \eqref{eq:hom-tensor}.
These in turn correspond to a collection, parameterized by $n \in \N$, of natural linear maps
\begin{equation}\label{eq:natural_linear_maps}
	\set[\big]{\triangle_i \colon \chains(\simplex^n) \to \chains(\simplex^n)^{\ot 2}}_{i \in \N}
\end{equation}
satisfying
\[
\bd \circ \, \triangle_i + \triangle_i \circ \bd =
(1+T) \triangle_{i-1}
\]
with $\triangle_{-1} = 0$.
Each natural transformation $\triangle_i$ in \eqref{eq:natural_linear_maps} is completely determined by the set of tensor chains
\[
\set[\big]{\triangle_i[n] \in \chains(\simplex^n)^{\ot 2}_{n+i}}_{n\in\N}
\]
where $[n]$ also denotes the identity $[n] \to [n]$, the top degree basis element.
We will often present a semi-simplicial \mbox{cup-$i$} construction $\triangle$ as the set $\set{\triangle_i[n]}_{i,n\in\N}$.

\medskip In the simplicial context, the simplex associated to a codegeneracy map $\sigma_j \colon [n] \to [n-1]$ is degenerate in $\simplex^{n-1}$ so it is $0$ in $\chains(\simplex^{n-1})$.
Therefore,
\[
\chains(\sigma_j)^{\ot 2} \triangle_i[n] = \big( \triangle_i \circ \chains(\sigma_j) \big)[n] = \triangle_i \big( \sigma_j \circ [n] \big) = \triangle_i(\sigma_j) = \triangle_i(0) = 0.
\]
which, together with the fact that codegeneracies and coface maps generate the simplex category, proves the following.

\begin{lemma}\label{l:simplicial_from_semisimplicial}
	A semi-simplicial \mbox{cup-$i$} construction $\triangle$ is simplicial if and only if $\chains(\sigma_j)^{\ot 2}\big(\triangle_i[n]\big) = 0$ for all $i,n\in\N$ and each codegeneracy map $\sigma_j \colon [n] \to [n-1]$.
\end{lemma}

\subsection{A dual description}\label{ss:equivalence_of_functors}

For $n,m \in \N$ let $\P_{n-m}^n$ be the set of subsets of $\set{0,\dots,n}$ whose cardinality is $n-m$.
Let $\cP(\simplex^n)$ be the chain complex defined by
\[
\cP(\simplex^n)_m = \Ftwo\set{\P_{n-m}^n}, \qquad
\bd U =
\begin{cases}
	\displaystyle\sum_{\mathclap{\bar{u} \notin U}} \bar{u}.U, & m>0, \\
	\hfil 0, & m=0,
\end{cases}
\]
where $\bar{u} \notin U$ stands for $\bar{u} \in \set{0,\dots,n} \setminus U$ and $\bar{u}.U$ for $\set{\bar{u}} \union U$.
For any coface $\delta_j \colon [n] \to [n+1]$ the chain map $\cP(\delta_j) \colon \cP(\simplex^n) \to \cP(\simplex^{n+1})$ is defined by
\[
\cP(\delta_j)(U) = \set{u_1 < \dots < u_{k-1} < j < u_k+1 < \dots < u_{n-m}+1}
\]
where $k$ is determined by the inequalities.
For any codegeneracy $\sigma_j \colon [n] \to [n-1]$ the chain map $\cP(\sigma_j) \colon \cP(\simplex^n) \to \cP(\simplex^{n-1})$ is defined by
\begin{equation}\label{eq:codegeneracies}
	\cP(\sigma_j)(U) =
	\begin{cases}
		\sigma_j\big(U \setminus \{j+1\}\big), & j+1 \in U, \\
		\hfil \sigma_j\big(U \setminus \{j\}\big), & j+1 \notin U \text{ and } j \in U, \\
		\hfil 0, & j+1 \notin U \text{ and } j \notin U,
	\end{cases}
\end{equation}
where $\sigma_j$ also denotes the induced map on subsets, i.e., $\sigma_j(S) = \set{\sigma_j(u) \mid u \in S}$.

\begin{lemma}\label{st:chains_functor_equivalence}
	The function
	\[
	\begin{tikzcd}[column sep=small,row sep=0]
		\Psi \colon &[-19] \P_{n-m}^n \rar & \simplex^n_m \\
		& U \rar[mapsto] & d_U [n]
	\end{tikzcd}
	\]
	induces a natural equivalence between the functors $\cP$ and $\chains$ with both domains: semi-simplicial and simplicial sets.
\end{lemma}

\begin{proof}
	Since any non-degenerate simplex in $\simplex^n$ is a face of the top dimensional one $[n]$, this assignment defines a bijection between the basis of $\cP(\simplex^n)_m$ and $\chains(\simplex^n)_m$.
	Let $U = \{u_1 < \dots < u_{n-m}\} \in \rP_{n-m}^n$.
	For $m=0$, both sides have zero boundary.
	For $m>0$, using that for $u \leq j$ we have $d_jd_u = d_ud_{j+1}$, we have
	\begin{align*}
		\bd d_U[n] &=
		\sum_{\mathclap{j \in \{0, \dots, m\}}}
		d_j\, d_{u_1} \dots\, d_{u_{n-m}}[n] \\ &=
		\sum_{\mathclap{\bar u \in \{0, \dots, n\} \setminus U}}
		d_{u_1} \dots\, d_{\bar{u}} \dots\, d_{u_{n-m}}[n] \\ &=
		\sum_{\mathclap{\bar u \in \{0, \dots, n\} \setminus U}}
		d_{\{\bar u\} \union U}[n] \\ &=
		d_{\bd U}[n].
	\end{align*}
	The compatibility with cofaces and codegeneracies follows directly from their definitions.
\end{proof}

\subsection{Cup-$i$ constructions, dually}\label{ss:axioms_revisited}

By combining \cref{ss:naturality} and \cref{ss:equivalence_of_functors}, we can see that the data of a semi-simplicial \mbox{cup-$i$} construction $\triangle$ is a collection, indexed by $i,n \in \N$, of elements
\begin{equation*}\label{eq:cup-i_dually}
	\triangle_i [n] =
	\sum_{\mathclap{\quad V \ot W \in \Lambda(i,n)}} \, V \ot W
\end{equation*}
in $\cP(\simplex^n)^{\ot 2}_{i+n}$
where $\Lambda(i,n)$ denotes the (possibly empty) set of basis elements occurring as summands of $\triangle_i[n]$.
For a chain $\zeta$ and a basis element $V \ot W$ of $\cP(\simplex^n)^{\ot 2}$ we write $\coeff{\zeta}{V \ot W} \in \F$ for the coefficient of $V \ot W$ in $\zeta$.

\begin{lemma}\label{l:dual_axioms}
	A semi-simplicial \mbox{cup-$i$} construction $\triangle$ is:
	\begin{enumerate}
		\item Irreducible if and only if $V \cap W = \emptyset$
		for every $i, n \in \N$ and $V \ot W \in \Lambda(i,n)$.
		\item Free if and only if $V \ot W \in \Lambda(i,n)$ implies $W \ot V \notin \Lambda(i,n)$,
		for every $i, n \in \N$ with $i \neq n$.
	\end{enumerate}
\end{lemma}

\begin{proof}
	There is $k \in V \cap W$ for one of the summands $V \ot W$ of $\triangle_i [n]$ if and only if the image of this summand under $\Psi^{\ot 2}$ is
	\[
	\underbrace{d_{V \setminus k} \overbrace{d_k[n]}^{y}}_{y^{(1)}}
	\ot
	\underbrace{d_{W \setminus k} \overbrace{d_k[n]}^{y}}_{y^{(2)}}
	\]
	with $(y^{(1)} \smallsmile_i y^{(2)})[n] \neq 0$.
	This proves the claim about irreducibility; the one involving freeness follows from direct inspection.
\end{proof}
	\section{Existence}

Using the reformulation established in the previous section, we present a simplicial cup-\(i\) construction satisfying all our axioms.
It plays a central role in the proof of our main theorem in \cref{s:proof}, and in \cref{ss:original} we prove that it coincides with Steenrod's original construction.
This construction appeared first in \cite{medina2025fast_sq} and was used in \cite{medina2022per_st} to compute persistent Steenrod squares of real-world data.

\subsection{Canonical construction}

For any $U \in \rP^n_{n-i}$ the \defn{index function} of $U$ is given by
\[
\begin{split}
	\ind_U \colon U &\to \F \\
	u_j &\mapsto u_j + j \mod 2.
\end{split}
\]
For any $U \in \rP^n_{n-i}$ and $\varepsilon \in \Ftwo \cong \{0,1\}$ we write $U^\varepsilon$ instead of $\ind_U^{-1}(\varepsilon)$.
For $0 \leq i \leq n$, define
\[
\canonical_i[n] =
\sum_{\mathclap{\quad U \in \rP^n_{n-i}}} \ U^0 \ot U^1.
\]
For $i>n$, set $\canonical_i[n]=0$.
By \cite[Theorem~10]{medina2025fast_sq}, these elements define a semi-simplicial \mbox{cup-$i$} construction.

\begin{theorem}\label{t:existence}
	The canonical \mbox{cup-$i$} construction is simplicial, non-zero, irreducible, and free.
\end{theorem}

\begin{proof}
	It is clearly non-zero and a straightforward application of \cref{l:dual_axioms} establishes its irreducibility and freeness.
	To see it is simplicial we consider \cref{eq:codegeneracies} and an arbitrary basis element $U^0 \ot U^1$ appearing as a summand in $\canonical_i[n]$.
	Let $U = U^0 \union U^1$.
	Based on \cref{l:simplicial_from_semisimplicial}, we need to check that for any $j \in \set{0,\dots,n-1}$ we have
	\begin{equation}\label{eq:checking_simplicial}
		\cP(\sigma_j)(U^0) \ot \cP(\sigma_j)(U^1) = 0.
	\end{equation}
	If $j \notin U$ or $j+1 \notin U$ either $\cP(\sigma_j)(U^0) = 0$ or $\cP(\sigma_j)(U^1) = 0$, so Identity \eqref{eq:checking_simplicial} holds.
	Let us assume $j,j+1 \in U$.
	Since $\ind_U(j) = \ind_U(j+1)$ then either $j,j+1 \in U^0$ or $j,j+1 \in U^1$.
	In both cases we also have \eqref{eq:checking_simplicial}.
\end{proof}

\subsection{Special cases}\label{ss:canonical_special_cases}

We close by recording two important instances of the canonical construction.
For any $n \in \N$,
\[
\canonical_n[n] = \emptyset \ot \emptyset,
\]
and, for $n > 0$,
\[
\canonical_{n-1} [n] \ = \
\sum_{\mathclap{\substack{\quad u \in \{0,\dots,n\} \\ u \ \mathrm{odd}}}} \,\{u\} \ot \emptyset \ + \
\sum_{\mathclap{\substack{\quad u \in \{0,\dots,n\} \\ u \ \mathrm{even}}}} \,\emptyset \ot \{u\}.
\]
	\section{Irreducibility}\label{s:irreducible}

In this section, we study the irreducibility condition.
It will be convenient to work in the \defn{augmented simplex category}, which contains an additional object \([-1] = \emptyset\) and a unique morphism from it to every \([n]\).
The natural equivalence between $\cP$ and $\chains$ of \cref{st:chains_functor_equivalence} extends to this category; we denote the extended functors by $\widehat{\cP}$ and $\widehat\chains$, respectively.
For every $n$, the only additional basis element of $\widehat{\cP}(\simplex^n)$ is the full subset $\set{0,\dots,n}$ in degree $-1$.
We freely identify the underlying graded vector spaces of $\widehat{\cP}(\simplex^n)$ and $\cP(\simplex^n)$ in non-negative degrees.

The main results of this section are an explicit description, in terms of a Hopf algebra structure carried by $\widehat{\cP}(\simplex^n)$, of the largest subcomplex of irreducible chains, and a proof that irreducibility together with semi-simpliciality implies simpliciality.

\subsection{A Koszul model}

Recall that an element $a$ of a commutative $\F$-algebra $A$ with $a^2 = 0$ makes $(A,\, a \cdot -)$ into a chain complex, referred to as the \defn{Koszul complex} of $a$.
Consider the exterior algebra over $\F$ generated by $x_0, \dots, x_n$, one generator for each vertex of $\simplex^n$, graded so that the monomial $x_U$ associated to $U \subseteq \set{0,\dots,n}$ has degree $n - \bars{U}$, with $x_\emptyset = 1$.
The assignment $U \mapsto x_U$ defines an isomorphism of complexes between $\widehat{\cP}(\simplex^n)$ and the Koszul complex of
\[
\omega = x_0 + \dots + x_n.
\]
We use it to regard $\widehat{\cP}(\simplex^n)$ as an algebra.
In what follows we use $U$ and $x_U$ interchangeably.
This Koszul complex factors as the tensor product of the complexes $\big(\F[x_u]/(x_u^2),\, x_u \cdot -\big)$, each of them acyclic.

\begin{remark*}\label{r:group_ring}
	Writing $t_u = 1 + x_u$ identifies the algebra $\widehat{\cP}(\simplex^n)$ with the group ring $\F[C_2^{\, n+1}]$ of an elementary abelian $2$-group, under which the Koszul element $\omega = \sum_u (1 + t_u)$ is the sum of the augmentation elements of its factors.
	The differentials of the complex $W$ of \cref{s:statement} are of the same nature, being multiplication by the augmentation element $1 + T$ of $\F[C_2]$.
\end{remark*}

\subsection{The partition coproduct}

The exterior algebra is a Hopf algebra, with coproduct the algebra map determined by declaring the generators primitive: $\Delta(x_u) = x_u \ot 1 + 1 \ot x_u$.
Explicitly,
\[
\Delta(x_U) = \sum x_V \ot x_W,
\]
the sum taken over all $V$ and $W$ satisfying $V \union W = U$ and $V \cap W = \emptyset$.
We refer to $\zeta_U = \Delta(x_U) \in \widehat{\cP}(\simplex^n)^{\ot 2}$ as the \defn{partition chain} of $U$.
Since $\Delta$ is cocommutative, partition chains are fixed by $T$.
If $U \subsetneq \set{0,\dots,n}$, then $\zeta_U$ is concentrated in non-negative bidegrees and agrees with its image in the ordinary complex $\cP(\simplex^n)^{\ot 2}$.
With respect to the grading above, $\Delta$ raises degree by $n$, so it defines a map $\Sigma^n \widehat{\cP}(\simplex^n) \to \widehat{\cP}(\simplex^n)^{\ot 2}$ from the $n$-fold suspension.
The next lemma states that this map is a chain map.

\begin{lemma}\label{c:boundary_partition_chain}
	For every $U \subseteq \set{0,\dots,n}$ we have
	\begin{equation}\label{eq:boundary_partition_chain}
		\bd\zeta_U = \sum_{\bar u \notin U} \zeta_{\bar u.U}.
	\end{equation}
\end{lemma}

\begin{proof}
	The differential of $\widehat{\cP}(\simplex^n)^{\ot 2}$ is multiplication by $\Delta(\omega) = \omega \ot 1 + 1 \ot \omega$.
	Since $\Delta$ is an algebra map,
	\[
	\bd \Delta(x_U) = \Delta(\omega)\, \Delta(x_U) = \Delta(\omega\, x_U) = \sum_{\bar u \notin U} \Delta(x_{\bar u.U}). \qedhere
	\]
\end{proof}

The proof exhibits $\bd$ as a coderivation of $\Delta$, reflecting the primitivity of $\omega$.
It is, however, not a derivation of the product: $\bd(x_V x_W) = \omega\, x_V x_W$, whereas $\bd(x_V)\, x_W + x_V\, \bd(x_W) = 2\, \omega\, x_V x_W = 0$.
Thus, $\widehat{\cP}(\simplex^n)$ is a differential graded coalgebra but not a differential graded Hopf algebra.

\begin{remark*}
	The partition coproduct does \emph{not} correspond to the group-ring coproduct $\F[C_2^{\, n+1}]$ as it makes each $x_u$ primitive rather than each $t_u$ grouplike.
\end{remark*}

\subsection{Irreducible chains}

A basis element $x_V \ot x_W \in \widehat{\cP}(\simplex^n)^{\ot 2}$ is said to be \defn{irreducible} if $V \cap W = \emptyset$ (equivalently, if the product $x_Vx_W$ is non-zero) and \defn{reducible} otherwise.
We denote by $\pired$ and $\piirred$ the projections onto the graded subspaces generated by reducible and irreducible basis elements, respectively.
Under the identification above, a cup-$i$ construction $\triangle$ is irreducible if and only if its defining elements $\triangle_i[n]$ are irreducible in this sense for all $i,n \in \N$.

\begin{remark}
	Irreducibility admits a geometric interpretation.
	Under the equivalence of \cref{st:chains_functor_equivalence}, a basis element $x_V \ot x_W$ corresponds to the pair of faces $\big(d_V[n],\, d_W[n]\big)$, and it is irreducible precisely when
	\[
	\dim d_V[n] + \dim d_W[n] - \dim \big(d_V[n] \cap d_W[n]\big) = n,
	\]
	that is, when the two faces intersect transversally in $\simplex^n$, in which case the product $x_V x_W$ corresponds to their intersection $d_{V \union W}[n]$.
	The multiplication of $\widehat{\cP}(\simplex^n)$ is thus the mod~$2$ intersection product of transverse faces, and the irreducibility axiom states that a \mbox{cup-$i$} construction is supported on transverse pairs.
\end{remark}

The reducible basis elements span the ideal of $\widehat{\cP}(\simplex^n)^{\ot 2}$ generated by the elements $x_k \ot x_k$.
Since the differential of $\widehat{\cP}(\simplex^n)^{\ot 2}$ is multiplication by $\Delta(\omega)$, every ideal is closed under it; in particular, the reducible subspace is a subcomplex.
We remark that this ideal is properly contained in the kernel of the multiplication map $\widehat{\cP}(\simplex^n)^{\ot 2} \to \widehat{\cP}(\simplex^n)$, which also contains irreducible chains such as $x_i \ot x_j + x_j \ot x_i$.
The irreducible subspace is not a subcomplex, and we now describe the largest subcomplex it contains.

\begin{theorem}\label{st:augmented_partition_subcomplex}
	Partition chains span the largest subcomplex contained in the irreducible subspace of $\widehat{\cP}(\simplex^n)^{\ot 2}$.
\end{theorem}

\begin{proof}
	Let $I$ be the irreducible subspace of $\widehat{\cP}(\simplex^n)^{\ot 2}$ and let $A$ be the span of the partition chains.
	The largest subcomplex contained in \(I\) is \(I \cap \bd^{-1}(I)\).
	The summands of every partition chain are irreducible, so $A \subseteq I$, and $A$ is a subcomplex by \cref{c:boundary_partition_chain}; hence \(A \subseteq I \cap \bd^{-1}(I)\).
	For the other inclusion,
	let $\zeta = \sum_\Lambda V \ot W$ be in $I \cap \bd^{-1}(I)$, where $\Lambda$ is a set of distinct irreducible basis elements.
	The reducible part of the boundary of $\zeta$ is
	\begin{equation}\label{eq:reducible_augmented_boundary}
		\pired(\bd\zeta) =
		\sum_\Lambda\left(
		\sum_{w \in W} w.V \ot W +
		\sum_{v \in V} V \ot v.W
		\right).
	\end{equation}
	Distinct terms within each inner family are distinct, while a term of the first family agrees with one of the second precisely when the corresponding partitions differ by a \defn{basic move}: the transfer of a single element between the factors.
	Hence \cref{eq:reducible_augmented_boundary} vanishes if and only if $\Lambda$ is closed under basic moves.
	Since basic moves preserve unions and connect any two partitions of the same set, this occurs exactly when $\Lambda$ is a union of complete families of partitions, that is, when $\zeta \in A$.
\end{proof}

We will apply this result to chains in the ordinary complex.
If $\zeta \in \cP(\simplex^n)^{\ot 2}$ is a homogeneous irreducible chain of degree at least $n$, then the augmented and ordinary boundaries of $\zeta$ have the same reducible part: their difference consists of terms having the full subset as one factor, the degree assumption forces the other factor to be $\emptyset$, and such terms are irreducible.
Consequently, if $\pired(\bd\zeta) = 0$ then $\zeta$ is a sum of partition chains, necessarily indexed by sets with at most $n$ elements.

\subsection{Irreducibility implies simpliciality}

The following theorem is not used in this paper but it is included to clarify the relationship between two axioms: irreducibility and naturality.

\begin{theorem}\label{t:irreducible_implies_simplicial}
	Every irreducible semi-simplicial \mbox{cup-$i$} construction is simplicial.
\end{theorem}

\begin{proof}
	Let $\triangle$ be an irreducible semi-simplicial \mbox{cup-$i$} construction.
	By \cref{l:simplicial_from_semisimplicial}, it suffices to show that
	\[
		\cP(\sigma_j)^{\ot 2}\triangle_i[n] = 0
	\]
	for every codegeneracy $\sigma_j \colon [n] \to [n-1]$.
	We derive this from the reducible part of the cup-$i$ identity.
	Fix $i,n \in \N$ with $n > 0$.
	We use the coefficient notation $\coeff{\zeta}{V \ot W}$ introduced in \cref{ss:axioms_revisited} and remind the reader that each $\triangle_i[n]$ has degree $n+i\geq n$.


	We first claim that
	\begin{equation}\label{eq:reducible_part}
		\pired\big(\bd\triangle_i[n]\big) = \triangle_i\bd\,[n].
	\end{equation}
	Indeed, $\triangle_i\bd\,[n] = \sum_{k=0}^{n} \cP(\delta_k)^{\ot 2}\triangle_i[n-1]$ and, since every subset in the image of $\cP(\delta_k)$ contains $k$, the basis elements in the image of $\cP(\delta_k)^{\ot 2}$ are reducible.
	On the other hand, $(1+T)\triangle_{i-1}[n]$ is a sum of irreducible basis elements.
	Applying $\pired$ to the identity $\bd\triangle_i[n] + \triangle_i\bd\,[n] = (1+T)\triangle_{i-1}[n]$ proves \eqref{eq:reducible_part}.

	Let us compare coefficients on both sides of \eqref{eq:reducible_part}.
	By \cref{eq:reducible_augmented_boundary}, every basis element appearing on the left-hand side of \cref{eq:reducible_part} is of the form $X \ot Y$ with $X \cap Y = \set{k}$ for exactly one $k$, and the same holds on the right-hand side, since $\cP(\delta_k)$ adjoins $k$ to both factors of an irreducible basis element.
	Fix such an $X \ot Y$.
	Its coefficient on the left-hand side is
	\[
	\coeff{\triangle_i[n]}{X \setminus k \ot Y} + \coeff{\triangle_i[n]}{X \ot Y \setminus k}.
	\]
	For the right-hand side, notice that $\cP(\delta_m)^{\ot 2}(A \ot B) = X \ot Y$ forces $m = k$, and that $\cP(\delta_k)$ defines a bijection from the subsets of $\set{0,\dots,n-1}$ to the subsets of $\set{0,\dots,n}$ containing $k$, with inverse $S \mapsto \sigma_k(S \setminus k)$, where $\sigma_k$ is the elementwise relabeling of \cref{eq:codegeneracies}.
	Therefore, \cref{eq:reducible_part} is equivalent to the identity
	\begin{equation}\label{eq:coefficient_comparison}
		\coeff{\triangle_i[n]}{X \setminus k \ot Y} + \coeff{\triangle_i[n]}{X \ot Y \setminus k} =
		\coeff{\triangle_i[n-1]}{\sigma_k(X \setminus k) \ot \sigma_k(Y \setminus k)}
	\end{equation}
	holding for every basis element $X \ot Y$ with $X \cap Y = \set{k}$.

	Let $\tau$ denote the exchange of $j$ and $j+1$, acting on subsets of $\set{0,\dots,n}$.
	Inspecting \cref{eq:codegeneracies}, we see that $\cP(\sigma_j)(\tau U) = \cP(\sigma_j)(U)$ for every subset $U$, and that an irreducible basis element $V \ot W$ satisfies $\cP(\sigma_j)^{\ot 2}(V \ot W) \neq 0$ if and only if one of $j,j+1$ belongs to $V$ and the other belongs to $W$.
	Such basis elements come in pairs $\set[\big]{V \ot W,\ \tau(V) \ot \tau(W)}$ of distinct elements sharing the same non-zero image, and elements of different pairs have different images.
	Consequently, since $\triangle_i[n]$ is irreducible, $\cP(\sigma_j)^{\ot 2}\triangle_i[n] = 0$ if and only if
	\begin{equation}\label{eq:swap}
		\coeff{\triangle_i[n]}{V \ot W} = \coeff{\triangle_i[n]}{\tau(V) \ot \tau(W)}
	\end{equation}
	for every such element.

	To establish \cref{eq:swap}, assume that $j \in V$ and $j+1 \in W$, the other case being symmetric.
	Since $W$ contains $j+1$ but not $j$, the set $W \union \set{j}$ is fixed by $\tau$ and equals $\tau(W) \union \set{j{+}1}$, so the basis elements $V \ot \big(W \union \set{j}\big)$ and $\tau(V) \ot \big(W \union \set{j}\big)$ have factors intersecting in $\set{j}$ and $\set{j+1}$, respectively.
	Applying \cref{eq:coefficient_comparison} to each of them gives
	\begin{align*}
		\coeff{\triangle_i[n]}{V \setminus j \ot W \union \set{j}} &+ \coeff{\triangle_i[n]}{V \ot W} \\ &=
		\coeff{\triangle_i[n-1]}{\sigma_j(V \setminus j) \ot \sigma_j(W)}, \\
		\coeff{\triangle_i[n]}{V \setminus j \ot W \union \set{j}} &+ \coeff{\triangle_i[n]}{\tau(V) \ot \tau(W)} \\ &=
		\coeff{\triangle_i[n-1]}{\sigma_{j+1}(V \setminus j) \ot \sigma_{j+1}\big(\tau(W)\big)}.
	\end{align*}
	The maps $\sigma_j$ and $\sigma_{j+1}$ agree on subsets containing neither $j$ nor $j+1$, so the first factors on the right agree, and so do the second, both being $\sigma_j\big(W \setminus \set{j{+}1}\big) \union \set{j}$.
	Adding the two identities therefore yields \cref{eq:swap}.
\end{proof}


\section{Freeness}\label{s:freeness}

In this section we record the way the freeness axiom will be used: for a free construction, the element $\triangle_i[n]$ is determined by $(1+T)\triangle_i[n]$ up to the transposition of its summands.

\begin{lemma}\label{l:freeness}
	Let $\triangle$ be a free semi-simplicial \mbox{cup-$i$} construction and let $i, n \in \N$ with $i \neq n$.
	Consider a set $\Lambda$ of basis elements of $\cP(\simplex^n)^{\ot 2}_{i+n}$, none of which is symmetric and no two of which are related by transposition.
	If
	\[
	(1+T)\, \triangle_i[n] = (1+T) \sum_{b \in \Lambda} b,
	\]
	then there is a unique function $\xi \colon \Lambda \to \set{0,1}$ such that
	\[
	\triangle_i[n] = \sum_{b \in \Lambda} T^{\xi(b)}\, b.
	\]
\end{lemma}

\begin{proof}
	Freeness, as expressed in \cref{l:dual_axioms}, states that $\triangle_i[n]$ does not contain both $b$ and $Tb$ as summands for any basis element $b$; in particular, no symmetric basis element is a summand of $\triangle_i[n]$.
	Hence, for any basis element $b$ with $b \neq Tb$, the identity
	\[
	\coeff{(1+T)\,\triangle_i[n]}{b} = \coeff{\triangle_i[n]}{b} + \coeff{\triangle_i[n]}{Tb}
	\]
	gives the following two implications:
	if $\coeff{(1+T)\,\triangle_i[n]}{b} = 0$ then neither $b$ nor $Tb$ is a summand of $\triangle_i[n]$, whereas if $\coeff{(1+T)\,\triangle_i[n]}{b} = 1$ then exactly one of them is.
	By hypothesis, $\coeff{(1+T)\,\triangle_i[n]}{b} = 1$ if and only if $b \in \set{c,\, T c}$ for a necessarily unique $c \in \Lambda$.
	Therefore, for each $c \in \Lambda$ exactly one of $c$ and $T c$ is a summand of $\triangle_i[n]$, and $\triangle_i[n]$ has no other summands.
	Defining $\xi(c) \in \set{0,1}$ accordingly establishes the claim.
\end{proof}

Taking $\Lambda = \emptyset$ gives the first corollary below.
For the second, notice that the basis elements $U^0 \ot U^1$ with $U \in \P^n_{n-i}$ satisfy the hypotheses of \cref{l:freeness}: they are parameterized by $U = U^0 \union U^1$; none is symmetric, since $U^0 = U^1$ forces $U = \emptyset$, which is excluded by $i \neq n$; and no two of them are related by transposition.

\begin{corollary}\label{c:free_zero}
	Let $\triangle$ be a free semi-simplicial \mbox{cup-$i$} construction and let $i, n \in \N$ with $i \neq n$.
	If $(1+T)\,\triangle_i[n] = 0$, then $\triangle_i[n] = 0$.
\end{corollary}

\begin{corollary}\label{c:free_xi}
	Let $\triangle$ be a free semi-simplicial \mbox{cup-$i$} construction and let $i, n \in \N$ with $i \neq n$.
	If
	\[
	(1+T)\,\triangle_i[n] = (1+T)\sum_{\mathclap{\quad U \in \P^n_{n-i}}} \, U^0 \ot U^1,
	\]
	then there is a function $\xi \colon \P^n_{n-i} \to \set{0,1}$ such that
	\[
	\triangle_i[n] = \sum_{\mathclap{\quad U \in \P^n_{n-i}}} \, U^\xi \ot U^\barxi,
	\]
	where $\barxi$ denotes the complement of $\xi$.
\end{corollary}


\section{Main proof}\label{s:proof}

In this section we present the proof of our main theorem: any semi-simplicial \mbox{cup-$i$} construction that is non-zero, irreducible and free is isomorphic to the canonical one.
For the rest of this section we let $\triangle$ be one such construction.

Our goal is to prove that for every $i\in\N$ there is $\varepsilon_i\in\set{0,1}$ such that
\[
\triangle_i[n]=T^{\varepsilon_i}\canonical_i[n]
\]
for every $n\in\N$.
We will do so in three stages \(i \geq n\), \(i = n - 1\), and \(i \leq n-2\).

\subsection{($i \geq n$)}\label{ss:i_geq_n}

\begin{proof}
	For $i > n$ the vector space $\cP(\simplex^n)^{\ot 2}_{n+i}$ is zero, which proves the claim.
	For $i = n$ the only non-zero element in the vector space $\cP(\simplex^n)^{\ot 2}_{2n}$ is $\emptyset \ot \emptyset = \canonical_n[n]$, which is fixed by $T$, so we are left with showing that $\triangle_n[n] \neq 0$ for every $n \in \N$.
	We will use without further mention that every symmetric basis element has coefficient $0$ in any element of the image of $(1+T)$.

	Let us first show that, for any $n > 0$, the elements $\triangle_n[n]$ and $\triangle_{n-1}[n-1]$ are either both zero or both equal to $\emptyset \ot \emptyset$.
	Assume, for a contradiction, that $\triangle_{n-1}[n-1] = \emptyset \ot \emptyset$ and $\triangle_n[n] = 0$.
	Since $\triangle_n[n-1] = 0$ by the $i > n$ case, we have
	\begin{align*}
	(1+T) \triangle_{n-1} [n] =
	\bd \triangle_{n} [n] + \triangle_{n} \bd \, [n] = 0,
	\end{align*}
	which, by \cref{c:free_zero}, implies $\triangle_{n-1}[n] = 0$.
	From this and $\triangle_{n-1}[n-1] = \emptyset \ot \emptyset$ we obtain
	\[
	(1+T)\triangle_{n-2}[n] =
	\bd \triangle_{n-1}[n] + \triangle_{n-1} \bd \, [n] =
	\sum_{u = 0}^{n} \{u\} \ot \{u\},
	\]
	a contradiction.
	Assume now that $\triangle_{n-1}[n-1] = 0$ and $\triangle_n[n] = \emptyset \ot \emptyset$.
	Then
	\begin{align*}
		(1+T) \triangle_{n-1} [n] &=
		\bd \triangle_{n} [n] + \triangle_{n} \bd \, [n] =
		\sum_{u=0}^n \set{u} \ot \emptyset + \emptyset \ot \set{u} \\ &=
		(1+T)\sum_{U \in \P_1^n} U^0 \ot U^1,
	\end{align*}
	so, by \cref{c:free_xi}, there is $\xi \colon \P_1^n \to \set{0,1}$ with $\triangle_{n-1}[n] = \sum_{U \in \P_1^n} U^\xi \ot U^\barxi$.
	Since $\triangle_{n-1}[n-1] = 0$ we have $\triangle_{n-1}\bd\,[n] = 0$ and therefore
	\[
	(1+T)\triangle_{n-2}[n] = \bd \triangle_{n-1}[n].
	\]
	For each $u \in \set{0, \dots, n}$, the symmetric basis element $\set{u} \ot \set{u}$ appears exactly once in $\bd \triangle_{n-1}[n]$: the summand of $\triangle_{n-1}[n]$ supported on $\set{u}$, namely $\set{u} \ot \emptyset$ or $\emptyset \ot \set{u}$, contributes it via $\bd\, \emptyset = \sum_{\bar u} \set{\bar u}$, and no other summand does.
	This is again a contradiction.

	Therefore, either $\triangle_n[n] = \emptyset \ot \emptyset$ for every $n \in \N$, or $\triangle_n[n] = 0$ for every $n \in \N$.
	In the latter case, an induction argument over $n-i$ shows that $\triangle_i [n] = 0$ for every $i, n \in \N$, i.e., that $\triangle$ is the zero \mbox{cup-$i$} construction, contradicting the non-zero axiom:
	for the induction step, consider
	\begin{align*}
	(1+T) \triangle_{i-1} [n] =
	\bd \triangle_{i} [n] + \triangle_{i} \bd\, [n] = 0,
	\end{align*}
	which, by \cref{c:free_zero}, implies $\triangle_{i-1} [n] = 0$.
\end{proof}

\begin{remark*}
	We mention that irreducibility was not used in this proof.
\end{remark*}

\subsection{($i = n-1$)}\label{ss:i_equals_n_minus_1}

\begin{proof}
	We will proceed by induction on $n$.
	For $n = 0$, we have that $\triangle_{-1}[0] = 0$ which trivially agrees with both $\canonical_{-1}[0]$ and $T\canonical_{-1}[0]$.
	Using the $i \geq n$ case (\cref{ss:i_geq_n}) we have
	\begin{align*}
		(1+T)\triangle_{n-1}[n] &=
		\bd\triangle_n[n] + \triangle_n\bd[n] \\ &=
		\sum_{u=0}^n \set{u} \ot \emptyset + \emptyset \ot \set{u} \\ &=
		(1+T)\sum_{U \in \P_1^n} U^0 \ot U^1.
	\end{align*}
	By \cref{c:free_xi}, there is a function $\xi \colon \P_1^n \to \set{0,1}$ such that
	\begin{equation}\label{eq:triangle_n-1_n}
		\triangle_{n-1}[n] = \sum_{U \in \P_1^n} U^\xi \ot U^\barxi
	\end{equation}
	and we must show that either $\xi \equiv 0$ or $\xi \equiv 1$.
	Let us rewrite \cref{eq:triangle_n-1_n} as
	\begin{align*}
		\triangle_{n-1}[n] &=
		\canonical_{n-1}[n] \ + \ \sum_{\mathclap{\substack{\set{u} \in \P_1^n \\ \xi \set{u} \neq 0}}} \set{u} \ot \emptyset + \emptyset \ot \set{u} \\ &=
		\canonical_{n-1}[n] \ + \ \sum_{\mathclap{\substack{\set{u} \in \P_1^n \\ \xi \set{u} \neq 0}}} \zeta_{\set{u}}
	\end{align*}
	where $\zeta_{\set{u}}$ is the partition chain of the singleton $\set{u}$.
	We claim that
	\begin{equation}\label{eq:bd_triangle}
		\bd\sum_{\mathclap{\substack{\set{u} \in \P_1^n \\ \xi \set{u} \neq 0}}} \ \zeta_{\set{u}} = 0.
	\end{equation}
	By the induction hypothesis, there is $\varepsilon \in \set{0,1}$ such that $\triangle_{n-2}[n-1] = T^\varepsilon \canonical_{n-2}[n-1]$.
	Consider the irreducible chain
	\[
	\delta = \triangle_{n-2}[n] + T^\varepsilon \canonical_{n-2}[n].
	\]
	Using \cref{ss:i_geq_n} to identify $\triangle_{n-1}\bd\,[n]$ with $\canonical_{n-1}\bd\,[n]$, we have
	\begin{align*}
		(1+T)\delta &=
		(1+T)\triangle_{n-2}[n] + (1+T)\canonical_{n-2}[n] \\ &=
		\big(\bd\triangle_{n-1}[n] + \triangle_{n-1}\bd\,[n]\big) +
		\big(\bd\canonical_{n-1}[n] + \canonical_{n-1}\bd\,[n]\big) \\ &=
		\bd\big(\triangle_{n-1}[n] + \canonical_{n-1}[n]\big) \\ &=
		\bd \sum_{\mathclap{\substack{\set{u} \in \P_1^n \\ \xi \set{u} \neq 0}}} \ \zeta_{\set{u}},
	\end{align*}
	so \cref{eq:bd_triangle} is equivalent to $(1+T)\delta = 0$.
	If $n = 1$ this holds trivially since $\delta = 0$.
	Assume $n \geq 2$.
	Since $T$ commutes with the maps induced by coface maps, the induction hypothesis gives $\triangle_{n-2}\bd\,[n] = T^\varepsilon \canonical_{n-2}\bd\,[n]$, and therefore
	\begin{align*}
		(1+T)\triangle_{n-3}[n] &=
		\bd\triangle_{n-2}[n] + \triangle_{n-2}\bd\,[n] \\ &=
		\bd\big(T^\varepsilon \canonical_{n-2}[n] + \delta\big) + T^\varepsilon \canonical_{n-2}\bd\,[n] \\ &=
		T^\varepsilon (1+T)\canonical_{n-3}[n] + \bd\delta \\ &=
		(1+T)\canonical_{n-3}[n] + \bd\delta.
	\end{align*}
	Since $(1+T)\triangle_{n-3}[n]$ and $(1+T)\canonical_{n-3}[n]$ are in the kernel of $\pired$, so is $\bd\delta$; that is, the irreducible chain $\delta$ is in the kernel of $(\pired \circ \bd)$.
	The chain $\delta$ has degree $2n-2\geq n$, so \cref{st:augmented_partition_subcomplex} applies.
	It follows that $\delta$ is a sum of partition chains and, since partition chains are fixed by $T$, we conclude that $(1+T)\delta = 0$, proving \cref{eq:bd_triangle}.

	We now show that \cref{eq:bd_triangle} implies either $\xi \equiv 0$ or $\xi \equiv 1$.
	By \cref{eq:boundary_partition_chain}, the left hand side of \cref{eq:bd_triangle} equals
	\[
	\sum_{\mathclap{\substack{\set{u} \in \P_1^n \\ \xi \set{u} \neq 0}}} \
	\sum_{\ \bar{u} \neq u} \zeta_{\set{u} \union \set{\bar u}},
	\]
	in which the partition chain of a two element set $\set{u, \bar u}$ appears twice, and hence cancels, if $\xi\set{u} = \xi\set{\bar u} = 1$, and appears exactly once if $\xi\set{u} \neq \xi\set{\bar u}$.
	Arguing by contradiction, let us assume the existence of $u$ and $\bar{u}$ with $\xi\set{u} = 1$ and $\xi\set{\bar u} = 0$.
	Then, the partition chain $\zeta_{\set{u, \bar u}}$ survives in the sum above and, since partition chains of distinct sets share no basis elements, \cref{eq:bd_triangle} is violated, a contradiction.
	Therefore, either $\triangle_{n-1}[n] = \canonical_{n-1}[n]$ or $\triangle_{n-1}[n] = T\canonical_{n-1}[n]$.
\end{proof}

\begin{remark*}\label{r:freeness_usage}
	The proof above uses freeness only through its dependence on \cref{ss:i_geq_n} and a single application of \cref{c:free_xi}, and the latter can be justified using irreducibility instead:
	in the relevant degree, the intersection of $\ker(1+T)$ with the subspace of irreducible chains is spanned by the partition chains $\zeta_{\set{u}}$, and adding these to $\triangle_{n-1}[n]$ only modifies the function $\xi$.
	The dependence on \cref{ss:i_geq_n} is, however, essential: the non-free construction $(1+T)\canonical$ of \cref{ss:non_free} is non-zero and irreducible, yet $(1+T)\canonical_n[n] = 0$ for every $n \in \N$.
\end{remark*}

\subsection{($i \leq n-2$)}\label{ss:proof}

\begin{proof}
	For each $i \in \N$, let $\varepsilon_i \in \set{0,1}$ be provided by \cref{ss:i_equals_n_minus_1} applied to $n = i+1$, so that $\triangle_i[i+1] = T^{\varepsilon_i}\canonical_i[i+1]$.
	We prove the statement with this choice of $\varepsilon_i$ using an induction argument over $k = n-i$.
	For $k \leq 0$ we have $\triangle_i [n] = \canonical_ i [n] = T^{\varepsilon_i}\canonical_ i [n]$ by \cref{ss:i_geq_n}, since these elements are fixed by $T$, and for $k = 1$ the claim holds by the choice of $\varepsilon_i$.
	Let $k > 1$, assume that the claim holds for all $i,n \in \N$ with $n-i < k$, and fix $i,n \in \N$ with $n-i = k$.
	Applying the induction hypothesis to $\triangle_{i+1}[n]$ and $\triangle_{i+1}[n-1]$ we have
	\begin{align*}
		(1+T)\triangle_i[n] &=
		\bd\triangle_{i+1}[n] + \triangle_{i+1}\bd[n] \\ &=
		T^{\varepsilon_{i+1}} \big( \bd\canonical_{i+1}[n] + \canonical_{i+1}\bd[n] \big) \\ &=
		T^{\varepsilon_{i+1}} (1+T)\,\canonical_i[n] \\ &=
		(1+T)\sum_{\mathclap{\quad U \in \P_{n-i}^n}} \, U^0 \ot U^1.
	\end{align*}
	\cref{c:free_xi} implies the existence of a function $\xi \colon \P_{n-i}^n \to \set{0,1}$ with
	\begin{equation}\label{eq:triangle_i_n}
		\triangle_i[n] = \sum_{\mathclap{\quad U \in \P_{n-i}^n}} \, U^{\xi} \ot U^{\barxi}.
	\end{equation}
	Consider the irreducible chain
	\[
	\delta = \triangle_i[n] + T^{\varepsilon_i}\canonical_i[n] =
	(1+T) \sum_{\quad\mathclap{\substack{U \in \P_{n-i}^n \\ \xi(U) \neq \varepsilon_i}}} U^0 \ot U^1,
	\]
	where the second identity follows from \cref{eq:triangle_i_n}: the summands with $\xi(U) = \varepsilon_i$ cancel, and the remaining ones contribute both $U^0 \ot U^1$ and $U^1 \ot U^0$.
	By the induction hypothesis $\triangle_i[n-1] = T^{\varepsilon_i}\canonical_i[n-1]$ and, since $T$ commutes with the maps induced by coface maps, $\triangle_{i}\bd[n] = T^{\varepsilon_i}\canonical_{i}\bd[n]$.
	Then
	\begin{equation}\label{eq:main proof a}
		\begin{split}
			(1+T)\triangle_{i-1}[n] &=
			\bd \triangle_{i}[n] \ +\ \triangle_{i}\bd[n] \\ &=
			\bd \big( T^{\varepsilon_i}\canonical_{i}[n] + \delta \big)
			\ +\ T^{\varepsilon_i}\canonical_{i}\bd[n] \\ &=
			T^{\varepsilon_i}(1+T)\canonical_{i-1}[n] \ +\ \bd\delta \\ &=
			(1+T)\canonical_{i-1}[n] \ +\ \bd\delta.
		\end{split}
	\end{equation}
	Since $(1+T)\triangle_{i-1}[n]$ and $(1+T)\canonical_{i-1}[n]$ are in the kernel of $\pired$, the irreducible chain $\delta$ is in the kernel of $(\pired \circ \bd)$.
	The chain $\delta$ has degree $n+i\geq n$, so \cref{st:augmented_partition_subcomplex} applies and shows that it is a sum of partition chains.
	If the partition chain of some $U \in \P^n_{n-i}$ appeared in this sum then, since partition chains of distinct sets share no basis elements, $\delta$ would contain all $2^{n-i}$ basis elements $V \ot W$ with $V \union W = U$ and $V \cap W = \emptyset$; but by its explicit form above, $\delta$ contains at most two of them, and $2^{n-i} > 2$ since $n - i = k > 1$.
	Therefore $\delta = 0$, that is, $\triangle_i[n] = T^{\varepsilon_i}\canonical_i[n]$, completing the induction.

	Finally, by the description of the automorphisms of $W$ given in \cref{s:statement}, the assignment $e_i \mapsto T^{\varepsilon_i} e_i$ defines an automorphism $\phi$ of $W$, inducing an isomorphism between the \mbox{cup-$i$} constructions $\triangle$ and $\canonical$.
	This completes the proof of the main theorem.
\end{proof}


\section{Independence of the axioms}\label{s:examples}

In this section we show that the axioms of our main theorem are independent: for each of them we exhibit a semi-simplicial \mbox{cup-$i$} construction satisfying the other two that is not isomorphic to the canonical one.
All examples are in fact simplicial.
Throughout this section we use the description of semi-simplicial \mbox{cup-$i$} constructions given in \cref{ss:axioms_revisited}.

\subsection{Removing non-zero}

The zero \mbox{cup-$i$} construction satisfies the irreducibility and freeness axioms vacuously, and it is not isomorphic to the canonical construction since an isomorphism of \mbox{cup-$i$} constructions preserves the property of being zero.
Therefore, the non-zero axiom cannot be removed from our main theorem.

\subsection{Removing free}\label{ss:non_free}

Consider the simplicial \mbox{cup-$i$} construction
\[
(1+T)\canonical = \set[\big]{\canonical_i[n] + T\canonical_i[n]}_{i,n}.
\]
It is a simplicial \mbox{cup-$i$} construction: naturality and simpliciality are preserved by $1+T$, while
\[
	\bd(1+T)\canonical_i + (1+T)\canonical_i\bd
	= (1+T)^2\canonical_{i-1}
	= (1+T)\big((1+T)\canonical_{i-1}\big) = 0.
\]
It is non-zero; indeed, by \cref{ss:canonical_special_cases}, for any $n > 0$
\[
(1+T)\canonical_{n-1}[n] = \sum_{u=0}^{n} \set{u} \ot \emptyset \ + \ \emptyset \ot \set{u}.
\]
It is irreducible since the canonical construction is (\cref{t:existence}) and transposition preserves the irreducibility of basis elements.
It is not free; for example, both $\set{1} \ot \emptyset$ and $\emptyset \ot \set{1}$ are summands of $(1+T)\canonical_0[1]$.
Finally, an isomorphism of \mbox{cup-$i$} constructions replaces each $\triangle_i$ by either $\triangle_i$ or $T\triangle_i$, so it preserves freeness.
Since the canonical construction is free and $(1+T)\canonical$ is not, they are not isomorphic.
Therefore, the freeness axiom cannot be removed from our main theorem.

\subsection{Removing irreducible}

\begin{lemma}\label{st:transposing}
	Let $\triangle$ be a semi-simplicial \mbox{cup-$i$} construction and fix $i_0, n_0 \in \N$ with $i_0 > 0$.
	For all non-negative integers $i$ and $n$ define
	\[
	\widetilde\triangle_i[n] =
	\begin{cases}
		T\triangle_{i_0}[n_0] & \text{if }i=i_0 \text{ and } n=n_0, \\
		\big(\triangle_{i_0-1} + \triangle_{i_0}\bd\big)[n_0+1] & \text{if }i=i_0-1 \text{ and } n=n_0+1, \\
		\big(\triangle_{i_0-1} + \bd\triangle_{i_0}\big)[n_0] & \text{if }i=i_0-1 \text{ and } n=n_0, \\
		\triangle_{i}[n] & \text{otherwise}.
	\end{cases}
	\]
	Then $\widetilde\triangle = \set[\big]{\widetilde\triangle_i[n]}_{i,n}$ is a semi-simplicial \mbox{cup-$i$} construction.
\end{lemma}

\begin{proof}
	We need to verify that for all $i$ and $n$ the identity
	\[
	\big(\bd \widetilde\triangle_{i} + \widetilde\triangle_{i}\bd\big)[n] =
	(1+T) \widetilde\triangle_{i-1}[n]
	\]
	holds, knowing it does so for $\triangle$.
	A consequence of this identity that we will use without further mention is the following
	\[
	(1+T)\triangle_i\bd = (1+T)\bd\triangle_i.
	\]
	We will split the proof into five cases
	\[
	i < i_0-1 \quad;\quad i = i_0-1 \quad;\quad i = i_0 \quad;\quad i = i_0+1 \quad;\quad i > i_0+1.
	\]
	For $i > i_0+1$ the identity on $\widetilde\triangle$ reduces to that for $\triangle$ for all values of $n$.
	For $i = i_0+1$ there is one value of $n$ to consider.
	If $n = n_0$ we have
	\begin{align*}
		\big(\bd \widetilde\triangle_{i_0+1} + \widetilde\triangle_{i_0+1} \bd\big) &=
		\big(\bd \triangle_{i_0+1} + \triangle_{i_0+1} \bd\big) \\ &=
		(1+T)\triangle_{i_0} \\ &=
		(1+T)T\triangle_{i_0} \\ &=
		(1+T)\widetilde\triangle_{i_0}.
	\end{align*}
	For $i = i_0$ there are two values of $n$ to consider.
	For $n = n_0+1$ we have
	\begin{align*}
		\big(\bd \widetilde\triangle_{i_0} + \widetilde\triangle_{i_0} \bd\big) &=
		\big(\bd \triangle_{i_0} + T\triangle_{i_0} \bd\big) \\ &=
		\big((1+T)\bd\triangle_{i_0} + T(\bd\triangle_{i_0} + \triangle_{i_0}\bd)\big) \\ &=
		(1+T)\big(\triangle_{i_0}\bd \,+\, \triangle_{i_0-1}\big) \\ &=
		(1+T)\widetilde\triangle_{i_0-1}
	\end{align*}
	and for $n = n_0$ we have
	\begin{align*}
		\big(\bd \widetilde\triangle_{i_0} + \widetilde\triangle_{i_0} \bd\big) &=
		\big(\bd T \triangle_{i_0} + \triangle_{i_0} \bd\big) \\ &=
		\big((1+T)\bd\triangle_{i_0} + \bd\triangle_{i_0} + \triangle_{i_0}\bd\big) \\ &=
		(1+T)\big(\bd\triangle_{i_0} + \triangle_{i_0-1}\big) \\ &=
		(1+T)\widetilde\triangle_{i_0-1}.
	\end{align*}
	For $i = i_0-1$ there are three cases to consider.
	For $n = n_0+2$ we have
	\begin{align*}
		\big(\bd\widetilde\triangle_{i_0-1} + \widetilde\triangle_{i_0-1} \bd\big) &=
		\big(\bd\triangle_{i_0-1} + (\triangle_{i_0-1} + \triangle_{i_0} \bd)\bd\big) \\ &=
		(1+T)\triangle_{i_0-2} \\ &=
		(1+T)\widetilde\triangle_{i_0-2}.
	\end{align*}
	For $n = n_0+1$ we have
	\begin{align*}
		\big(\bd\widetilde\triangle_{i_0-1} + \widetilde\triangle_{i_0-1} \bd\big) &=
		\big(\bd(\triangle_{i_0-1} + \triangle_{i_0} \bd) + (\triangle_{i_0-1} + \bd\triangle_{i_0})\bd \big) \\ &=
		(1+T)\triangle_{i_0-2} \\ &=
		(1+T)\widetilde\triangle_{i_0-2}.
	\end{align*}
	For $n = n_0$ we have
	\begin{align*}
		\big(\bd\widetilde\triangle_{i_0-1} + \widetilde\triangle_{i_0-1} \bd\big) &=
		\big(\bd(\triangle_{i_0-1} + \bd\triangle_{i_0}) + \triangle_{i_0-1}\bd \big) \\ &=
		(1+T)\triangle_{i_0-2} \\ &=
		(1+T)\widetilde\triangle_{i_0-2}.
	\end{align*}
	If $i<i_0-1$ the identity on $\widetilde\triangle$ reduces to that for $\triangle$ for all values of $n$.
\end{proof}

Let us consider the canonical \mbox{cup-$i$} construction $\canonical$ and the integers $i_0=1$ and $n_0=3$.
The \mbox{cup-$i$} construction $\widetilde\triangle$ defined in \cref{st:transposing} is both non-zero and free but not irreducible as we will inspect.
If $i \neq i_0-1$ or $n \not\in \set{n_0,n_0+1}$ then $\widetilde\triangle_i[n]$ is equal to $\canonical_i[n]$ or $T\canonical_i[n]$.
With the assistance of \href{https://comch.readthedocs.io/en/latest/}{\texttt{ComCH}} we compute the two cases where the axioms could be broken:
\[
\widetilde\triangle_0[3] = (\canonical_0 + \partial\canonical_1)[3] =
\]
\noindent
{\ttfamily
	(0)(23) + (2)(23) + (3)(23) + ()(123) + (01)(3) + (13)(3) +\\
	(1)(13) + (012)() + (12)(1) + (12)(2) + (0)(03) + (3)(03) +\\
	(02)(0) + (2)(02) + (0)(01) + (1)(01).}
\[
\widetilde\triangle_0[4] = (\canonical_0 + \canonical_1\partial)[4] =
\]
\noindent
{\ttfamily
	(1234)() + (234)(0) + (34)(01) + (4)(012) + ()(0123) + (0)(034) +\\
	(02)(04) + (023)(0) + (0)(014) + (03)(01) + (0)(012) + (1)(134) +\\
	(12)(14) + (123)(1) + (1)(014) + (13)(01) + (1)(012) + (2)(234) +\\
	(12)(24) + (123)(2) + (2)(024) + (23)(02) + (2)(012) + (3)(234) +\\
	(13)(34) + (123)(3) + (3)(034) + (23)(03) + (3)(013) + (4)(234) +\\
	(14)(34) + (124)(4) + (4)(034) + (24)(04) + (4)(014).}

\medskip\noindent
The construction is non-zero because $\widetilde\triangle_0[0]=\canonical_0[0]\neq0$, and it is not irreducible since, for example, $(2)(23)$ is a summand of $\widetilde\triangle_0[3]$.
It is free because neither exceptional chain displayed above contains a summand together with its transposition, while every other component is canonical or its transposition.
It is also a simplicial \mbox{cup-$i$} construction.
By \cref{l:simplicial_from_semisimplicial}, we need each modified component to be annihilated by $\cP(\sigma_j)^{\ot 2}$ for every codegeneracy $\sigma_j$.
For $T\canonical_{i_0}[n_0]$ and $\bd\canonical_{i_0}[n_0]$ this holds because $\canonical$ is simplicial and $\cP(\sigma_j)^{\ot 2}$ commutes with $T$ and $\bd$.
For $\canonical_{i_0}\bd\,[n_0{+}1] = \sum_{m} \cP(\delta_m)^{\ot 2}\, \canonical_{i_0}[n_0]$ it follows from the cosimplicial identities: the summands with $m \notin \set{j, j{+}1}$ are annihilated because $\canonical$ is simplicial, whereas the two summands with $m \in \set{j, j{+}1}$ are equal and cancel since $\sigma_j \delta_j = \sigma_j \delta_{j+1} = \id$.
Since transposition preserves the irreducibility of basis elements, an isomorphism of \mbox{cup-$i$} constructions preserves irreducibility; the canonical construction is irreducible and $\widetilde\triangle$ is not, so they are not isomorphic.
Therefore, the irreducibility axiom cannot be removed from our main theorem.


\section{Other formulas}\label{s:others}

Several constructions of \mbox{cup-$i$} products appear in the literature.
In this section we show that they all agree with the canonical one, and therefore with each other.
We proceed in two ways.
For Steenrod's original construction (\cref{ss:original}) and for the one obtained from the Alexander--Whitney--Eilenberg--Zilber contraction (\cref{ss:real}) we verify the axioms of our main theorem directly and then determine which representative of the isomorphism class each of them is.
The remaining constructions are of operadic origin, and for these we do not repeat that analysis: we identify each of them with Steenrod's construction, using the description of the latter in terms of the Alexander--Whitney diagonal and the join (\cref{ss:join}), first for the surjection operad (\cref{ss:surjections}) and then, through table reduction, for the Barratt--Eccles operad (\cref{ss:barratt_eccles}).
The agreement of these constructions seems to appear in print for the first time here.

The constructions reviewed below are defined with integer coefficients, where their formulas acquire signs.
Since our comparison takes place with $\Ftwo$-coefficients, we do not engage with these signs, referring to the cited references for the general case.

\subsection{Original construction}\label{ss:original}

We now review, using the perspective developed in \cref{s:reformulation}, \defn{Steenrod's \mbox{cup-$i$} construction} \cite[p.293]{steenrod1947products}.
If $i > n$ then $\triangle_i^\rS [n] = 0$, otherwise it is given by the sum over all ordered sequences of integers
\[
0 \leq p_1 < \dots < p_{i+1} \leq n
\]
of the basis element
\begin{equation}\label{e:original i odd}
\begin{split}
[ 0, \dots, &{p_1} ] \ast [ {p_2}, \dots, {p_3} ] \ast \dots \ast [ {p_{i+1}}, \dots, n ] \\
\ot \ [ &{p_1}, \dots, {p_2} ] \ast \dots \ast [ {p_{i}}, \dots, {p_{i+1}} ]
\end{split}
\end{equation}
if $i$ is odd, and of
\begin{equation}\label{e:original i even}
\begin{split}
[ 0, \dots, &{p_1} ] \ast [ {p_2}, \dots, {p_3} ] \ast \dots \ast [ {p_{i}}, \dots, {p_{i+1}} ] \\
\ot \ [ &{p_1}, \dots, {p_2} ] \ast [ {p_3}, \dots, {p_4} ] \ast \dots \ast [ {p_{i+1}}, \dots, n ]
\end{split}
\end{equation}
if $i$ is even, where $\ast$ denotes the join of simplices:
\[
[{p_{k-1}}, \dots, {p_{k}} ] \ast [ {p_{k+1}}, \dots, p_{k+2}] = [{p_{k-1}}, \dots, p_k, p_{k+1}, \dots, p_{k+2}].
\]

\begin{theorem*}\label{t:steenrod cup-i}
	Steenrod's \mbox{cup-$i$} construction agrees with the canonical one.
\end{theorem*}

\begin{proof}
	We will first use our main result to prove that they are isomorphic.
	We can then conclude their equality by inspecting that $\triangle^\rS_i [i+1] \neq T \canonical_i [i+1]$.

	Steenrod's simplicial \mbox{cup-$i$} construction is \emph{non-zero} since $\triangle^\rS_0\big([0]\big) = [0] \ot [0] \neq 0$.
	It is \emph{irreducible} since for each basis element in $\triangle^\rS_i [n]$ with $i \leq n$, all integers $\{0, \dots, n\}$ appear in at least one of the tensor factors.
	To prove it is \emph{free} let us assume $i$ is even with $i < n$.
	The case where $i$ is odd is done analogously.
	If it is not free, then there exist two, not necessarily distinct, sequences
	\begin{align*}
	\begin{split}
	0 &= p_0 \leq p_1 < \dots < p_{i+1} \leq p_{i+2} = n \\
	0 &= q_0 \leq q_1 \,< \dots < q_{i+1} \leq q_{i+2} = n
	\end{split}
	\end{align*}
	such that
	\[
	\begin{split}
	&[ {p_0}, \dots, {p_1} ] \ast [ {p_2}, \dots, {p_3} ] \ast \dots \ast [ {p_{i}}, \dots, {p_{i+1}} ]\ = \\
	&[ {q_1}, \dots, {q_2} ] \ast [ {q_3}, \dots, {q_4} ] \ast \dots \ast [ {q_{i+1}}, \dots, {q_{i+2}} ]
	\end{split}
	\]
	and
	\[
	\begin{split}
	&[ {q_0}, \dots, {q_1} ] \ast [ {q_2}, \dots, {q_3} ] \ast \dots \ast [ {q_{i}}, \dots, {q_{i+1}} ]\ = \\
	&[ {p_1}, \dots, {p_2} ] \ast [ {p_3}, \dots, {p_4} ] \ast \dots \ast [ {p_{i+1}}, \dots, {p_{i+2}} ].
	\end{split}
	\]
	We will prove that $p_{r+1} = q_{r+1} = r$ for $0 \leq r \leq i$.
	Granting this, the largest vertices appearing in the two sides of the second identity are $q_{i+1} = i$ and $p_{i+2} = n$, respectively, giving the contradiction $i = n$.
	The base case $r = 0$ holds since the smallest vertices of the two sides of the first identity are $p_0 = 0$ and $q_1$, and those of the second are $q_0 = 0$ and $p_1$.
	For the induction step, assume $p_s = q_s = s - 1$ for $1 \leq s \leq r$ with $r \leq i$, and let $r$ be even.
	The vertices of the right hand side of the first identity smaller than $q_{r+1}$ are those of $[q_1, \dots, q_2] \ast \dots \ast [q_{r-1}, \dots, q_r]$, namely $0, \dots, r-1$.
	Since the left hand side contains the interval $[p_r, \dots, p_{r+1}]$ and $p_{r+1} > p_r = r - 1$, the vertex $r$ belongs to it, forcing $q_{r+1} = r$.
	The second identity gives $p_{r+1} = r$ in the same way, and for $r$ odd the same argument applies with the roles of the two identities exchanged.

	Our main theorem then shows that for each $i \in \N$ either $\triangle^\rS_i = \canonical_i$ or $\triangle^\rS_i = T \canonical_i$.
	Consider the element $U = \{0\} \in \rP_{1}^{i+1}$ giving rise to the summand $U^1 \ot U^0 = \set{0} \ot \emptyset$ in $T \canonical_i [i+1]$.
	Applying the isomorphism $\cP(\simplex^n)^{\ot 2} \to \chains(\simplex^n)^{\ot 2}$ we obtain the basis element $[1,\dots,i+1] \otimes [0,\dots,i+1]$ which is not a summand of $\triangle^\rS_i [i+1]$.
	This concludes the proof.
\end{proof}

\subsection{AW--EZ contraction}\label{ss:real}

In work by Real \cite{real1996computability}, further developed by Gonz\'alez\-/D\'iaz\--Real \cite{gonzalez-diaz1999steenrod, gonzalez2003computation, gonzalez-diaz2005cocyclic}, alternative formulas defining a simplicial \mbox{cup-$i$} construction were introduced using the Alexander--Whitney and Eilenberg--Zilber linear natural transformations
\[
\AW \colon \chains(\simplex^n \times \simplex^n)
\rightleftarrows
\chains(\simplex^n) \ot \chains(\simplex^n) \,: \EZ
\]
and an explicit natural chain homotopy $\mathrm{SHI}$ between their non-trivial composition $\EZ \AW$ and the identity.
The natural linear transformations defining \defn{Real's simplicial \mbox{cup-$i$} construction} are
\[
\triangle^{\mathrm{R}}_i = \AW (T \, \mathrm{SHI})^i.
\]

In \cite[Theorem~3.1]{gonzalez-diaz1999steenrod} these authors unravelled the above definition in terms of face maps, which we now review using the perspective developed in \cref{s:reformulation}.
If $i > n$ then $\triangle_i^{\mathrm{R}} [n] = 0$, otherwise it is given by
\begin{multline}\label{eq:real left}
	\triangle^{\mathrm{R}}_i [n] = \!
	\sum_{j_i=i}^{n} \ \sum_{j_{i-1}=i-1}^{j_i-1} \dots \sum_{j_0=0}^{j_1-1} \\
	\{j_0+1,\dots,j_1-1\} \union \{j_2+1,\dots,j_3-1\} \union \dots \union \{j_i+1,\dots,n\} \\ \ot \,
	\{0,\dots,j_0-1\} \union \{j_1+1,\dots,j_2-1\} \union \dots \union \{j_{i-1}+1,\dots,j_i-1\}
\end{multline}
if $i$ is even and by
\begin{multline}\label{eq:real right}
	\triangle^\rR_i [n] = \!
	\sum_{j_i=i}^{n} \ \sum_{j_{i-1}=i-1}^{j_i-1} \dots \sum_{j_0=0}^{j_1-1} \\
	\{j_0+1,\dots,j_1-1\} \union \{j_2+1,\dots,j_3-1\} \union \dots \union \{j_{i-1}+1,\dots,j_i-1\} \\ \ot \,
	\{0,\dots,j_0-1\} \union \{j_1+1,\dots,j_2-1\} \union \dots \union \{j_i+1,\dots,n\}
\end{multline}
if $i$ is odd.\footnote{To compare these formulas to those appearing in \cite[Theorem~3.1]{gonzalez-diaz1999steenrod} replace $i$ and $n$ respectively by $n$ and $m$.}

\begin{theorem*}
	Real's \mbox{cup-$i$} construction agrees with the canonical one.
\end{theorem*}

\begin{proof}
	We will use our main theorem to prove that they are isomorphic.
	We can then conclude their equality by inspecting that $\triangle^\rR_i [i+1] \neq T \canonical_i [i+1]$.
	Let us start by describing the indexing set of the sum in recursive terms.
	It consists of all integers $0 \leq j_0 < \dots < j_i \leq n$ with
	\begin{equation}\label{eq:indexing set recursively}
		j_i \in \set{i,\dots,n} \text{ and }
		j_k \in \set{k,\dots,j_{k+1}-1} \text{ for } k \in \set{0,\dots,i-1}.
	\end{equation}

	We can see that $\triangle^\rR$ is non-zero since, if $i=n=0$, the indexing set of the sum consists of a single $j_0 = 0$ and $\triangle^\rR_0 [0] = \emptyset \ot \emptyset$, using \cref{eq:real left}.

	We can see that $\triangle^\rR$ is irreducible since for each summand $V \ot W$, equal to
	\begin{align*}
		&\{j_0+1,\dots,j_1-1\} \union \{j_2+1,\dots,j_3-1\} \union \dots \union \{j_i+1,\dots,n\} \\ \ot \,
		&\{0,\dots,j_0-1\} \union \{j_1+1,\dots,j_2-1\} \union \dots \union \{j_{i-1}+1,\dots,j_i-1\}
	\end{align*}
	when $i$ is even, or
	\begin{align*}
		&\{j_0+1,\dots,j_1-1\} \union \{j_2+1,\dots,j_3-1\} \union \dots \union \{j_{i-1}+1,\dots,j_i-1\} \\ \ot \,
		&\{0,\dots,j_0-1\} \union \{j_1+1,\dots,j_2-1\} \union \dots \union \{j_i+1,\dots,n\}
	\end{align*}
	when $i$ is odd, we have $V \cap W = \emptyset$ given that $0 \leq j_0 < \dots < j_i \leq n$.

	To prove that $\triangle^\rR$ is free, let us focus on the case when $i$ is an even number, the odd case being analyzed similarly.
	Fix integers $0 \leq i \leq n$ and let us assume that $V \ot W = W' \ot V'$ for some pair of summands.
	After suppressing the union symbol and the commas this identity is explicitly given by
	\begin{align*}
		&\{j_0+1 \dots j_1-1\} \{j_2+1 \dots j_3-1\} \dots \{j_i+1 \dots n\} \\=
		&\{0 \dots j'_0-1\} \{j'_1+1 \dots j'_2-1\} \dots \{j'_{i-1}+1 \dots j'_i-1\}
	\end{align*}
	and
	\begin{align*}
		&\{j'_0+1 \dots j'_1-1\} \{j'_2+1 \dots j'_3-1\} \dots \{j'_i+1 \dots n\} \\=
		&\{0 \dots j_0-1\} \{j_1+1 \dots j_2-1\} \dots \{j_{i-1}+1 \dots j_i-1\}.
	\end{align*}
	Using the first of these identities we have that: if $j_i+1 \leq n$ then $n \in \{j_i+1 \dots n\}$ so $j'_i-1 \geq n$, which is impossible since $j'_i \in \set{i,\dots,n}$.
	Therefore, $j_i = n$.
	Using the other identity the same argument shows that $j'_i = n$.
	Let us consider this to be the base case of an induction argument on $k = i,i-1,\dots,0$, with hypothesis $j_k = j'_k = n-i+k$.
	Let us consider the induction step from $k$ to $k-1$.
	In this case,
	\[
	\set{j_{k+\ell}+1,\dots,j_{k+\ell+1}-1} = \emptyset = \set{j'_{k+\ell}+1,\dots,j'_{k+\ell+1}-1}
	\]
	for all $\ell \in \set{0,\dots,i-k-1}$.
	Let us assume $k$ is even.
	The first identity becomes
	\begin{align*}
		&\{j_0+1 \dots j_1-1\} \{j_2+1 \dots j_3-1\} \dots \{j_{k-2}+1 \dots j_{k-1}-1\} \\=
		&\{0 \dots j'_0-1\} \{j'_1+1 \dots j'_2-1\} \dots \{j'_{k-1}+1 \dots j'_k-1\}.
	\end{align*}
	If $j'_{k-1}+1 \leq j'_k-1$, then $j'_k - 1 = (n-i+k) - 1$ belongs to the right hand side and hence to the left hand side.
	The largest possible element of the left hand side is $j_{k-1}-1$, so this would imply $j_{k-1}-1\geq j'_k-1=j_k-1$, contradicting $j_{k-1}<j_k$.
	Therefore the last interval on the right is empty, and \cref{eq:indexing set recursively} gives $j'_{k-1}=j'_k-1=n-i+k-1$.
	A similar argument using the other identity shows that $j_{k-1} = n-i+k-1$.
	If $k$ is odd, the same argument applies with the two identities interchanged.
	Finally, using that $j_0 = j'_0 = n-i$, the first identity reduces to
	\[
	\emptyset = \set{0,\dots, n-i-1}
	\]
	implying $n-i \leq 0$, and, since $n-i \geq 0$, we conclude $n=i$.

	It remains to be shown that $\triangle^\rR$ is not only isomorphic to the canonical \mbox{cup-$i$} construction $\canonical$, but equal.
	To do so, notice that for every non-negative integer $i$ the indexing element with $j_k = k+1$ for each $k \in \{0,\dots,i\}$ is the unique one producing the summand $\emptyset \ot \{0\}$ in $\triangle_i^\rR[i+1]$.
	This is not a summand of $T \canonical_i[i+1]$.
	Therefore, $\triangle_i^\rR \neq T \canonical_i$ for every $i \geq 0$, so $\triangle^\rR = \canonical$ as claimed.
\end{proof}


\subsection{Diagonal and join}\label{ss:join}

The \defn{Alexander--Whitney diagonal} is the natural transformation
\[
\triangle_{\AWd} \colon \chains(X) \to \chains(X)^{\ot 2}
\]
sending an $n$-simplex $x$ to $\sum_{j=0}^{n} d_{j+1} \dotsm d_n(x) \ot d_0 \dotsm d_{j-1}(x)$; it is the natural transformation determined, as in \cref{ss:naturality}, by the elements $\canonical_0[n]$.
The \defn{join} is the natural transformation
\[
\ast \colon \chains(\simplex^n)^{\ot 2} \to \chains(\simplex^n)
\]
of degree $1$ defined on basis elements by
\[
\ast \big([v_0, \dots, v_p] \ot [w_0, \dots, w_q]\big) =
\begin{cases}
	[u_0, \dots, u_{p+q+1}] &
	\begin{aligned}
		&\text{if } \set{v_0, \dots, v_p} \cap {}\\[-2pt]
		&\hspace{17pt}\set{w_0, \dots, w_q} = \emptyset,
	\end{aligned}\\
	\hfil 0 & \text{if not},
\end{cases}
\]
where $u_0 < \dots < u_{p+q+1}$ are the vertices involved.

The following recursive description was recorded in \cite{medina2023dennis}; we include a direct verification from Steenrod's formulas.

\begin{lemma}
	Steenrod's construction admits the recursive description
	\begin{equation}\label{eq:prop cup-i}
		\triangle^\rS_0 = \triangle_{\AWd},
		\qquad
		\triangle^\rS_i = (\ast \ot \id) \circ (\id \ot T\triangle^\rS_{i-1}) \circ \triangle_{\AWd}.
	\end{equation}
\end{lemma}

\begin{proof}
	We compare this recursion inductively with the formulas of \cref{ss:original}.
	For $i=0$, \cref{e:original i even} gives
	\[
		\triangle^\rS_0[n]
		=
		\sum_{p_1=0}^{n}
		[0,\dots,p_1] \ot [p_1,\dots,n]
		=
		\triangle_{\AWd}[n].
	\]
	Now fix $i>0$ and assume the formula in \cref{ss:original} has been recovered for $\triangle^\rS_{i-1}$.
	The Alexander--Whitney diagonal first chooses a vertex $p_1$ and produces
	\[
		[0,\dots,p_1] \ot [p_1,\dots,n].
	\]
	By the inductive hypothesis, applying $\triangle^\rS_{i-1}$ to the second factor gives a sum indexed by
	\[
		p_1 \leq p_2 < \dots < p_{i+1} \leq n.
	\]
	After transposing its two outputs, the first one is joined to $[0,\dots,p_1]$.
	The factor being joined has smallest vertex $p_2$.
	Consequently, the join vanishes when $p_2=p_1$ and is non-zero when $p_1<p_2$.
	Thus the surviving terms are indexed precisely by
	\[
		0 \leq p_1 < \dots < p_{i+1} \leq n.
	\]

	If $i$ is even, then $i-1$ is odd, and the surviving tensor is
	\[
	\begin{split}
		&[0,\dots,p_1] \ast [p_2,\dots,p_3] \ast \dots \ast [p_i,\dots,p_{i+1}] \\
		&\qquad \ot
		[p_1,\dots,p_2] \ast [p_3,\dots,p_4] \ast \dots \ast [p_{i+1},\dots,n].
	\end{split}
	\]
	This is the summand of \cref{e:original i even} indexed by the same sequence.
	If $i$ is odd, then $i-1$ is even, and the surviving tensor is
	\[
	\begin{split}
		&[0,\dots,p_1] \ast [p_2,\dots,p_3] \ast \dots \ast [p_{i+1},\dots,n] \\
		&\qquad \ot
		[p_1,\dots,p_2] \ast [p_3,\dots,p_4] \ast \dots \ast [p_i,\dots,p_{i+1}],
	\end{split}
	\]
	the corresponding summand of \cref{e:original i odd}.
	This proves the recursion.
\end{proof}

Steenrod's \mbox{cup-$i$} products are therefore composites of the Alexander--Whitney diagonal, the join, and the transposition of tensor factors.
As proved in \cite{medina2020prop1}, arbitrary such composites assemble into a natural $E_\infty$-coalgebra structure on simplicial chains, a structure induced from a cellular analogue \cite{medina2021prop2}.
The operads considered next act through it.

\subsection{Surjections}\label{ss:surjections}

Set
\[
\triangle_{\AWd}^{[1]} = \id, \quad
\triangle_{\AWd}^{[m+1]} = \big(\triangle_{\AWd}^{[m]} \ot \id\big) \circ \triangle_{\AWd},
\qquad
\ast^{[1]} = \id, \quad
\ast^{[m+1]} = \ast \circ \big(\ast^{[m]} \ot \id\big),
\]
with the coassociativity of $\triangle_{\AWd}$ and the associativity of $\ast$ making the chosen parenthesizations immaterial.
Let
\[
s \colon \set{1, \dots, \ell} \twoheadrightarrow \set{1, \dots, r}
\]
be a surjection, represented by the sequence $\big(s(1), \dots, s(\ell)\big)$.
In degree $\ell-r$, the basis of $\cX(r)$ consists of the \defn{nondegenerate surjections}, those satisfying $s(j) \neq s(j+1)$ for every $j$; surjections with consecutive repetitions represent zero.
Its action on simplicial chains is the natural transformation
\begin{equation}\label{e:surjection action}
	\Big(\ast^{[\bars{s^{-1}(1)}]} \ot \dotsb \ot \ast^{[\bars{s^{-1}(r)}]}\Big)
	\circ \pi_s \circ \triangle_{\AWd}^{[\ell]},
\end{equation}
where $\pi_s$ permutes the tensor factors so as to group them according to the fibers of $s$, preserving their order within each fiber.
These operations define the standard action of the McClure--Smith surjection operad $\cX$ \cite{mcclure2003multivariable}, realized inside the structure of \cref{ss:join} in \cite[Appendix~1]{medina2020prop1}.
The analogous description holds integrally with the appropriate signs, which we suppress since this paper works over \(\Ftwo\).

Consider the \defn{alternating surjections}
\[
(1,2), \qquad (1,2,1), \qquad (1,2,1,2), \qquad \dots
\]
Let $a_i^\cX$ be the alternating surjection with $i+2$ entries and let $A_i$ be its action on simplicial chains.
The action of $a_0^\cX = (1,2)$ is the Alexander--Whitney diagonal, so $A_0 = \triangle_{\AWd}$.
For $i>0$, the two fibers of $a_i^\cX$ alternate along the sequence; grouping them in \cref{e:surjection action}, with the last factor transposed, gives
\[
A_i = (\ast \ot \id) \circ (\id \ot T A_{i-1}) \circ \triangle_{\AWd}.
\]
Comparing this identity with \cref{eq:prop cup-i} and proceeding by induction shows that $A_i = \triangle^\rS_i$ for every $i$.
The alternating surjections therefore define Steenrod's \mbox{cup-$i$} construction.

\subsection{Barratt--Eccles}\label{ss:barratt_eccles}

The \defn{Barratt--Eccles operad} $\cE$ is given by $\cE(r) = \chains(E\Sym_r)$, with operad structure induced by the partial composition of permutations \cite{berger2004combinatorial}.
Berger and Fresse constructed a natural $\cE$-coalgebra structure on simplicial chains which factors through the surjection operad along their table reduction morphism $\TR \colon \cE \to \cX$.
For each $i \in \N$, let
\[
a_i^\cE = (\id, T, \id, \dots, T^i)
\]
be the alternating sequence with $i+1$ entries; these begin
\[
(\id), \qquad (\id, T), \qquad (\id, T, \id), \qquad \dots
\]
The table reduction morphism satisfies $\TR(a_i^\cE) = a_i^\cX$, so the alternating sequences also define Steenrod's \mbox{cup-$i$} construction.

We refer to \href{https://comch.readthedocs.io/en/latest/}{\texttt{ComCH}} for a computer algebra system implementing the surjection and Barratt--Eccles operads.

	\appendix
	\sloppy
	\printbibliography

@article{chen2012symmetry,
  author = {Chen, Xie and Gu, Zheng-Cheng and Liu, Zheng-Xin and Wen, Xiao-Gang},
  title = {Symmetry protected topological orders in interacting bosonic systems},
  journal = {Science},
  volume = {338},
  number = {6114},
  pages = {1604--1606},
  year = {2012},
  doi = {10.1126/science.1227224},
  eprint = {1301.0861},
  eprinttype = {arXiv},
  eprintclass = {cond-mat.str-el}
}

@article{chen2013symmetry,
  author = {Chen, Xie and Gu, Zheng-Cheng and Liu, Zheng-Xin and Wen, Xiao-Gang},
  title = {Symmetry protected topological orders and the group cohomology of their symmetry group},
  journal = {Physical Review B},
  volume = {87},
  number = {15},
  pages = {155114},
  year = {2013},
  doi = {10.1103/PhysRevB.87.155114},
  eprint = {1106.4772},
  eprinttype = {arXiv},
  eprintclass = {cond-mat.str-el}
}

@article{feng2026anyonic,
  author = {Feng, Yitao and Xue, Hanyu and Li, Yuyang and Cheng, Meng and Kobayashi, Ryohei and Hsin, Po-Shen and Chen, Yu-An},
  title = {Anyonic membranes and {P}ontryagin statistics},
  journal = {Physical Review Letters},
  volume = {136},
  number = {8},
  pages = {086601},
  year = {2026},
  doi = {10.1103/4jww-6b6t},
  eprint = {2509.14314},
  eprinttype = {arXiv},
  eprintclass = {cond-mat.str-el}
}

@article{gu2014symmetry,
  author = {Gu, Zheng-Cheng and Wen, Xiao-Gang},
  title = {Symmetry-protected topological orders for interacting fermions: Fermionic topological nonlinear {$\sigma$} models and a special group supercohomology theory},
  journal = {Physical Review B},
  volume = {90},
  number = {11},
  pages = {115141},
  year = {2014},
  doi = {10.1103/PhysRevB.90.115141},
  eprint = {1201.2648},
  eprinttype = {arXiv},
  eprintclass = {cond-mat.str-el}
}

@article{kapustin2015fermionic_cobordisms,
  author = {Kapustin, Anton and Thorngren, Ryan and Turzillo, Alex and Wang, Zitao},
  title = {Fermionic symmetry protected topological phases and cobordisms},
  journal = {Journal of High Energy Physics},
  volume = {2015},
  number = {12},
  pages = {052},
  year = {2015},
  doi = {10.1007/JHEP12(2015)052},
  eprint = {1406.7329},
  eprinttype = {arXiv},
  eprintclass = {cond-mat.str-el}
}

@article{meng2018classification,
	title = {Classification of symmetry-protected phases for interacting fermions in two dimensions},
	author = {Cheng, Meng and Bi, Zhen and You, Yi-Zhuang and Gu, Zheng-Cheng},
	journal = {Phys. Rev. B},
	volume = {97},
	issue = {20},
	pages = {205109},
	numpages = {9},
	year = {2018},
	publisher = {American Physical Society},
	doi = {10.1103/PhysRevB.97.205109},
	url = {https://link.aps.org/doi/10.1103/PhysRevB.97.205109}
}

@article{wang2020construction,
	title = {Construction and Classification of Symmetry-Protected Topological Phases in Interacting Fermion Systems},
	author = {Wang, Qing-Rui and Gu, Zheng-Cheng},
	journal = {Phys. Rev. X},
	volume = {10},
	issue = {3},
	pages = {031055},
	numpages = {64},
	year = {2020},
	publisher = {American Physical Society},
	doi = {10.1103/PhysRevX.10.031055},
	url = {https://link.aps.org/doi/10.1103/PhysRevX.10.031055}
}

@article{medina2024connected,
	author = {Cantero-Mor{\'a}n, Federico and Medina-Mardones, Anibal M.},
	title = {Connected power operations and simplicial Poincar\'e duality},
	journal = {arXiv e-prints},
	year = {2024},
	archivePrefix = {arXiv},
	eprint = {2402.00826},
	doi = {10.48550/arXiv.2402.00826},
	url = {https://arxiv.org/abs/2402.00826},
}

@article {medina2025odd_cartan,
	AUTHOR = {Cantero-Mor\'{a}n, Federico and Medina-Mardones, Anibal M.},
	TITLE = {An effective proof of the {C}artan formula: odd primes},
	JOURNAL = {Homology Homotopy Appl.},
	FJOURNAL = {Homology, Homotopy and Applications},
	VOLUME = {27},
	YEAR = {2025},
	NUMBER = {1},
	PAGES = {207--234},
	ISSN = {1532-0073},
	MRCLASS = {55S10 (55S05 55S12)},
	MRNUMBER = {4883682},
	DOI = {10.4310/hha.2025.v27.n1.a12},
	URL = {https://doi.org/10.4310/hha.2025.v27.n1.a12},
}

@article {medina2024multisimplicial,
	AUTHOR = {Medina-Mardones, Anibal M. and Pizzi, Andrea and Salvatore,
	Paolo},
	TITLE = {Multisimplicial chains and configuration spaces},
	JOURNAL = {J. Homotopy Relat. Struct.},
	FJOURNAL = {Journal of Homotopy and Related Structures},
	VOLUME = {19},
	YEAR = {2024},
	NUMBER = {2},
	PAGES = {275--296},
	ISSN = {2193-8407},
	MRCLASS = {55R80 (18N40 18N50 18N70 55U15)},
	MRNUMBER = {4746152},
	DOI = {10.1007/s40062-024-00344-7},
	URL = {https://doi.org/10.1007/s40062-024-00344-7},
}

@article {medina2023dennis,
	AUTHOR = {Medina-Mardones, Anibal M.},
	TITLE = {The diagonal of cellular spaces and effective algebro-homotopical constructions},
	JOURNAL = {EMS Surv. Math. Sci.},
	FJOURNAL = {EMS Surveys in Mathematical Sciences},
	VOLUME = {10},
	YEAR = {2023},
	NUMBER = {2},
	PAGES = {223--241},
	ISSN = {2308-2151},
	MRCLASS = {55 (18)},
	MRNUMBER = {4667420},
	DOI = {10.4171/emss/71},
	URL = {https://doi.org/10.4171/emss/71},
}

@article{medina2025fast_sq,
	AUTHOR = {{Medina-Mardones}, Anibal M.},
	TITLE = {New formulas for cup-{$i$} products and fast computation of {S}teenrod squares},
	JOURNAL = {Comput. Geom.},
	FJOURNAL = {Computational Geometry. Theory and Applications},
	VOLUME = {109},
	YEAR = {2023},
	ISSN = {0925-7721},
	MRNUMBER = {4473678},
	DOI = {10.1016/j.comgeo.2022.101921},
	URL = {https://doi.org/10.1016/j.comgeo.2022.101921},
}

@article{medina2022cube_einfty,
	author = {{Kaufmann}, Ralph M. and {Medina-Mardones}, Anibal M.},
	title = "{A combinatorial ${E}_\infty$-algebra structure on cubical cochains and the Cartan--Serre map}",
	journal = {Cahiers Topologie G\'{e}om. Diff\'{e}rentielle Cat\'{e}g.},
	fjournal = {Cahiers de Topologie et G{\'e}om{\'e}trie Diff{\'e}rentielle Cat{\'e}goriques},
	volume = {63},
	number = {4},
	pages = {387--424},
	year = {2022},
	url = {http://cahierstgdc.com/wp-content/uploads/2022/10/KAUFMANN-MEDINA-LXIII-4.pdf},
}

@article{medina2022per_st,
	author = {{Lupo}, Umberto and {Medina-Mardones}, Anibal M. and {Tauzin}, Guillaume},
	title = "{Persistence Steenrod modules}",
	fjournal = {Journal of Applied and Computational Topology},
	journal = {J. Appl. Comput. Topol.},
	volume={6},
	number={4},
	pages={475--502},
	year = {2022},
	doi = {10.1007/s41468-022-00093-7},
	URL = {https://doi.org/10.1007/s41468-022-00093-7},
	publisher = {Springer}
}

@ARTICLE{medina2021comch,
	author = {Anibal M. {Medina-Mardones}},
	title = {{A computer algebra system for the study of commutativity up to coherent homotopies}},
	volume = {14},
	journal = {Advanced Studies: Euro-Tbilisi Mathematical Journal},
	number = {4},
	publisher = {Tbilisi Centre for Mathematical Sciences},
	pages = {147 -- 157},
	year = {2021},
	url = {https://projecteuclid.org/journals/advanced-studies-euro-tbilisi-mathematical-journal/volume-14/issue-4/A-computer-algebra-system-for-the-study-of-commutativity-up/10.3251/asetmj/1932200819.full}
}

@ARTICLE{medina2021prop2,
	AUTHOR = {{Medina-Mardones}, Anibal M.},
	TITLE = {A finitely presented {$E_\infty$}-prop {II}: cellular context},
	JOURNAL = {High. Struct.},
	FJOURNAL = {Higher Structures},
	VOLUME = {5},
	YEAR = {2021},
	NUMBER = {1},
	PAGES = {69--186},
	URL = {https://higher-structures.math.cas.cz/api/files/issues/Vol5Iss1/Medina-Mardones-2}
}

@article {medina2021may_st,
	AUTHOR = {Kaufmann, Ralph M. and Medina-Mardones, Anibal M.},
	TITLE = {Cochain level {M}ay-{S}teenrod operations},
	JOURNAL = {Forum Math.},
	FJOURNAL = {Forum Mathematicum},
	VOLUME = {33},
	YEAR = {2021},
	NUMBER = {6},
	PAGES = {1507--1526},
	ISSN = {0933-7741},
	MRCLASS = {55S05 (55S10 55S12 55S15 55U05 55U15)},
	MRNUMBER = {4333989},
	DOI = {10.1515/forum-2020-0296},
	URL = {https://doi.org/10.1515/forum-2020-0296},
}

@article {medina2021adem,
	AUTHOR = {Brumfiel, Greg and {Medina-Mardones}, Anibal and Morgan, John},
	TITLE = {A cochain level proof of {A}dem relations in the mod 2 {S}teenrod algebra},
	JOURNAL = {J. Homotopy Relat. Struct.},
	FJOURNAL = {Journal of Homotopy and Related Structures},
	VOLUME = {16},
	YEAR = {2021},
	NUMBER = {4},
	PAGES = {517--562},
	ISSN = {2193-8407},
	MRCLASS = {55S10},
	MRNUMBER = {4343073},
	DOI = {10.1007/s40062-021-00287-3},
	URL = {https://doi.org/10.1007/s40062-021-00287-3},
}

@ARTICLE{medina2020prop1,
	AUTHOR = {{Medina-Mardones}, Anibal M.},
	TITLE = {A finitely presented {$E_\infty$}-prop {I}: algebraic context},
	JOURNAL = {High. Struct.},
	FJOURNAL = {Higher Structures},
	VOLUME = {4},
	YEAR = {2020},
	NUMBER = {2},
	PAGES = {1--21},
	MRCLASS = {55P48 (18M85 18N50)},
	MRNUMBER = {4133162},
	url = {https://journals.mq.edu.au/api/files/issues/Vol4Iss2/Medina-Mardones}
}

@ARTICLE{medina2020cartan,
	AUTHOR = {{Medina-Mardones}, Anibal M.},
	TITLE = {An effective proof of the {C}artan formula: the even prime},
	JOURNAL = {J. Pure Appl. Algebra},
	FJOURNAL = {Journal of Pure and Applied Algebra},
	VOLUME = {224},
	YEAR = {2020},
	NUMBER = {12},
	PAGES = {106444, 18},
	ISSN = {0022-4049},
	MRCLASS = {55S10 (55S05 55S12)},
	MRNUMBER = {4102178},
	DOI = {10.1016/j.jpaa.2020.106444},
	URL = {https://doi.org/10.1016/j.jpaa.2020.106444},
}

@article {mcclure2003multivariable,
	AUTHOR = {McClure, James E. and Smith, Jeffrey H.},
	TITLE = {Multivariable cochain operations and little {$n$}-cubes},
	JOURNAL = {J. Amer. Math. Soc.},
	FJOURNAL = {Journal of the American Mathematical Society},
	VOLUME = {16},
	YEAR = {2003},
	NUMBER = {3},
	PAGES = {681--704},
	ISSN = {0894-0347},
	MRCLASS = {55P48 (18D50)},
	MRNUMBER = {1969208},
	MRREVIEWER = {Benoit Fresse},
	DOI = {10.1090/S0894-0347-03-00419-3},
	URL = {https://doi.org/10.1090/S0894-0347-03-00419-3},
}

@article {berger2004combinatorial,
	AUTHOR = {Berger, Clemens and Fresse, Benoit},
	TITLE = {Combinatorial operad actions on cochains},
	JOURNAL = {Math. Proc. Cambridge Philos. Soc.},
	FJOURNAL = {Mathematical Proceedings of the Cambridge Philosophical
	Society},
	VOLUME = {137},
	YEAR = {2004},
	NUMBER = {1},
	PAGES = {135--174},
	ISSN = {0305-0041},
	MRCLASS = {18D50 (16E45 55P48)},
	MRNUMBER = {2075046},
	MRREVIEWER = {David Chataur},
	DOI = {10.1017/S0305004103007138},
	URL = {https://doi.org/10.1017/S0305004103007138},
}

@article {steenrod1947products,
	AUTHOR = {Steenrod, N. E.},
	TITLE = {Products of cocycles and extensions of mappings},
	JOURNAL = {Ann. of Math. (2)},
	FJOURNAL = {Annals of Mathematics. Second Series},
	VOLUME = {48},
	YEAR = {1947},
	PAGES = {290--320},
	ISSN = {0003-486X},
	MRCLASS = {56.0X},
	MRNUMBER = {22071},
	MRREVIEWER = {B. Eckmann},
	DOI = {10.2307/1969172},
	URL = {https://doi.org/10.2307/1969172},
}

@article {real1996computability,
	AUTHOR = {Real, Pedro},
	TITLE = {On the computability of the {S}teenrod squares},
	JOURNAL = {Ann. Univ. Ferrara Sez. VII (N.S.)},
	FJOURNAL = {Annali dell'Universit\`a di Ferrara. Nuova Serie. Sezione VII.
	Scienze Matematiche},
	VOLUME = {42},
	YEAR = {1996},
	PAGES = {57--63 (1998)},
	ISSN = {0430-3202},
	MRCLASS = {55S05 (55S10)},
	MRNUMBER = {1622630},
	MRREVIEWER = {N. J. Kuhn},
	URL = {https://doi.org/10.1007/BF02955020}
}

@incollection {gonzalez-diaz1999steenrod,
	AUTHOR = {Gonz\'{a}lez-D\'{i}az, Roc\'{i}o and Real, Pedro},
	TITLE = {A combinatorial method for computing {S}teenrod squares},
	NOTE = {Effective methods in algebraic geometry (Saint-Malo, 1998)},
	JOURNAL = {J. Pure Appl. Algebra},
	FJOURNAL = {Journal of Pure and Applied Algebra},
	VOLUME = {139},
	YEAR = {1999},
	NUMBER = {1-3},
	PAGES = {89--108},
	ISSN = {0022-4049},
	MRCLASS = {55N45 (55S10 55U10 55U15)},
	MRNUMBER = {1700539},
	MRREVIEWER = {Kalathoor Varadarajan},
	DOI = {10.1016/S0022-4049(99)00006-7},
	URL = {https://doi.org/10.1016/S0022-4049(99)00006-7},
}

@article{gonzalez2003computation,
	title={Computation of cohomology operations of finite simplicial complexes},
	author={Gonzalez-Diaz, Rocio and Real, Pedro and others},
	journal={Homology, Homotopy and Applications},
	volume={5},
	number={2},
	pages={83--93},
	year={2003},
	publisher={International Press of Boston},
	url={https://dx.doi.org/10.4310/HHA.2003.v5.n2.a4}
}

@article {gonzalez-diaz2005cocyclic,
	AUTHOR = {{Gonzalez-Diaz}, Rocio and Real, Pedro},
	TITLE = {{HPT} and cocyclic operations},
	JOURNAL = {Homology Homotopy Appl.},
	FJOURNAL = {Homology, Homotopy and Applications},
	VOLUME = {7},
	YEAR = {2005},
	NUMBER = {2},
	PAGES = {95--108},
	ISSN = {1532-0081},
	MRCLASS = {55S10 (05E99)},
	MRNUMBER = {2156309},
	MRREVIEWER = {Adriana Ciampella},
	URL = {http://projecteuclid.org/euclid.hha/1139839376},
}

@incollection {pilarczyk2016cubical,
	AUTHOR = {Kr\v{c}\'{a}l, Marek and Pilarczyk, Pawe{\l} },
	TITLE = {Computation of cubical {S}teenrod squares},
	BOOKTITLE = {Computational topology in image context},
	SERIES = {Lecture Notes in Comput. Sci.},
	VOLUME = {9667},
	PAGES = {140--151},
	PUBLISHER = {Springer},
	YEAR = {2016},
	MRCLASS = {55N45 (55U15 68W30)},
	MRNUMBER = {3533883},
	DOI = {10.1007/978-3-319-39441-1\_13},
	URL = {https://doi.org/10.1007/978-3-319-39441-1_13},
}

@article {gaiotto2016spin,
	author = {Gaiotto, Davide and Kapustin, Anton},
	title = {Spin TQFTs and fermionic phases of matter},
	journal = {International Journal of Modern Physics A},
	volume = {31},
	number = {28n29},
	pages = {1645044},
	year = {2016},
	doi = {10.1142/S0217751X16450445},
	URL = {https://doi.org/10.1142/S0217751X16450445},
	eprint = {https://doi.org/10.1142/S0217751X16450445},
}

@article {brumfiel2016pontrjagin,
	author = {{Brumfiel}, Greg and {Morgan}, John},
	title = "{The {P}ontrjagin Dual of 3-Dimensional Spin Bordism}",
	journal = {arXiv e-prints},
	year = {2016},
	archivePrefix = {arXiv},
	eprint = {1612.02860},
}

@article {kapustin2017fermionic,
	AUTHOR = {Kapustin, Anton and Thorngren, Ryan},
	TITLE = {Fermionic {SPT} phases in higher dimensions and bosonization},
	JOURNAL = {J. High Energy Phys.},
	FJOURNAL = {Journal of High Energy Physics},
	YEAR = {2017},
	NUMBER = {10},
	PAGES = {080, front matter+48},
	ISSN = {1126-6708},
	MRCLASS = {81T30},
	MRNUMBER = {3731133},
	DOI = {10.1007/jhep10(2017)080},
	URL = {https://doi.org/10.1007/jhep10(2017)080},
}

@article {brumfiel2018pontrjagin,
	author = {{Brumfiel}, Greg and {Morgan}, John},
	title = "{The {P}ontrjagin Dual of 4-Dimensional Spin Bordism}",
	journal = {arXiv e-prints},
	year = {2018},
	archivePrefix = {arXiv},
	eprint = {1803.08147},
}

@article{barkeshli2022classification,
	title={Classification of (2+1)D invertible fermionic topological phases with symmetry},
	author={Barkeshli, Maissam and Chen, Yu-An and Hsin, Po-Shen and Manjunath, Naren},
	journal={Physical Review B},
	volume={105},
	number={23},
	pages={235143},
	year={2022},
	publisher={American Physical Society},
	doi={10.1103/PhysRevB.105.235143},
	url={https://doi.org/10.1103/PhysRevB.105.235143}
}
	\enlargethispage{3\baselineskip}
	\todos
\end{document}